\documentclass[a4wide,11pt]{scrartcl}

\usepackage[USenglish]{babel}
\usepackage[T1]{fontenc}
\usepackage[utf8]{inputenc}

\usepackage{graphicx,subfigure,tikz,pgfplots,epstopdf}
\usepackage{amsmath,amssymb,amsfonts,amsthm,array,geometry,stmaryrd,booktabs,paralist,todonotes,mathtools,pdflscape,siunitx,placeins}

\usepackage[colorlinks=true,linkcolor=blue,citecolor=green,urlcolor=blue,bookmarks]{hyperref}

\numberwithin{equation}{section}

\newcommand{\R}{\mathbb{R}}
\newcommand{\dd}{\textup{d}}

\newtheorem{thm}{Theorem}[section]

\newtheorem{rem}[thm]{Remark}

\newenvironment{prf}{\textbf{Proof:}} {\hspace*{\fill} $\square$ \newline}

\usepackage{url}
\usepackage{dsfont}
\usepackage{nicefrac}
\usepackage{xspace}
\usepackage{framed}

\usepackage{cite}

\newcommand{\p}[2]{\ensuremath{\frac{\partial #1}{\partial #2 }}}
\newcommand{\diag}{\mathop{\mathrm{diag}}}

\begin{document}
    \title{On the relation between approaches for boundary feedback control of hyperbolic systems}
    \author{Michael Herty\footnotemark[1] \and  Ferdinand Thein\footnotemark[1]\;\;\footnotemark[2]}
    \date{}%\today}
    \maketitle
    \begin{abstract}
        Stabilization of partial differential equations is a topic of utmost importance in mathematics as well as in engineering sciences.
        Concerning one dimensional problems there exists a well developed theory.
        Due to numerous important applications the interest in boundary feedback control of multi-dimensional hyperbolic systems is increasing.
        In the present work we want to discuss the relation between some of the most recent results available in the literature. The key result of the present work is to show that the type of system discussed in \cite{Yang2024} identifies a particular class which falls into the framework presented in \cite{Herty2023a}.
    \end{abstract}
    \renewcommand{\thefootnote}{\fnsymbol{footnote}}
    \footnotetext[1]{IGPM RWTH Aachen, Templergraben 55, D-52056 Aachen, Germany.\\
    \href{mailto:herty@igpm.rwth-aachen.de}{\textit{herty@igpm.rwth-aachen.de}}}
    \footnotetext[2]{present address: Johannes Gutenberg-Universit\"at, Staudingerweg 9, D-55128 Mainz, Germany.\\
    \href{mailto:fthein@uni-mainz.de}{\textit{fthein@uni-mainz.de}}}
    \renewcommand{\thefootnote}{\arabic{footnote}}
    \section{Introduction}
    The stabilization of spatially one--dimensional systems of hyperbolic balance laws is a vivid subject attracting research interest in the mathematical as well as in the engineering community
    and we refer to the monographs \cite{Bastin2016,MR2302744,MR2655971,MR2412038} for further references. The mentioned references also provide a comprehensive overview on related controllability problems.
    A particular focus has been put on problems modeled by the barotropic Euler equations and the shallow water equations which in one space dimension form a 2$\times$2 hyperbolic system to model the temporal
    and spatial evolution of fluid flows including flows on networks.
    Analytical results concerning the boundary control of such systems have been studied in several articles, cf.\ \cite{G1,G2,G3,G4} for gas flow and for water flow we refer to \cite{W1,W2,W3,W4,W5,W6}.
    One key aspect in the analysis is the Lyapunov function which is introduced as a weighted $L^2$ (or $H^s$) norm and which allows to estimate deviations from steady states, see e.g. \cite{Bastin2016}.
    Under rather general {dissipative} conditions the exponential decay of the Lyapunov function has been established for various problem formulations and we exemplary refer to \cite{L1,L2,L4,L5}.
    For a study on comparisons to other stability concepts we mention \cite{L7}.
    Stability with respect to a higher $H^s$-norm $(s\geq 2)$ gives stability of the nonlinear system \cite{L5,Bastin2016}.
    Without aiming at completeness, we mention that recently the results have been extended to also deal with, e.g.,\ input-to-state stability \cite{MR2899713},
    numerical methods have been discussed in \cite{MR3956429,MR3031137,MR3648349} and for results concerning nonlocal hyperbolic partial differential equations (PDEs) see \cite{MR4172728}. A further widely used approach is the backstepping procedure introduced in \cite{Kanellakopoulos1991} and \cite{Kokotovic1992}. Due to the systematic procedure it gained interest and we for example refer to \cite{Ezal2000,Balogh2002,Smyshlyaev2005} for applications of the backstepping design to different problems. Furthermore in \cite{Vazquez2011} backstepping is applied to a two-by-two one dimensional hyperbolic system and very recently in \cite{Demir2024} the procedure is used to control a model for neuron growth. In a very recent paper backstepping is used together with linear matrix inequalities for the dynamic control of a system of interest, see \cite{Ibarra2024}. For further reading we refer to the aforementioned literature and the references therein. Despite its presence in the literature the further treatment of this method is beyond the scope of the present work.
    \par
    However, to the best of our knowledge the presented results are limited to the spatially one--dimensional case and multi--dimensional applications are ubiquitous.
    Thus there is a demand for a multi--dimensional extension of the available theory, yet results in the literature are rare due to the inherent difficulties and mostly focus on particular cases.
    Based on an application in metal forming processes, see \cite{Bambach2022,Herty2023}, we extended results to multi--dimensional hyperbolic balance laws which are simultaneously diagonalizable.
    In \cite{Herty2023a} an ansatz for symmetric hyperbolic systems is presented. It relies on the feasibility of an associated linear matrix inequality (LMI).
    A specific system in two dimensions is discussed in \cite{Dia2013} where a control problem for the shallow water equations is studied.
    There the authors take advantage of the structure of the system and show that the energy is non-increasing upon imposing suited boundary conditions.
    Just very recently another approach was proposed where the boundary stabilization for multi--dimensional systems is studied using a different Lyapunov function, see \cite{Yang2024}.
    There a stabilizing approach is introduced for a class of multi--dimensional systems satisfying the \emph{structural stability condition} (SSC).
    Concerning results on systems satisfying the SSC we exemplary refer to \cite{Yong1999,Yong2001,Yong2008,Herty2016,Yong2019} for further reading. As will be shown in this work the class of SSC systems is a particular representative for which the LMI condition holds.
    An additional recent result on multi--dimensional hyperbolic scalar conservation laws was presented in \cite{Serre2022}.
    However, the goal of the mentioned paper is different from the one presented here.
    Here we focus on the relation between the works \cite{Herty2023a} and \cite{Yang2024}. 
    We want to provide a comparative study which aims at giving a sound basis for further progress on this research subject.\\
    The outline of this paper is as follows.
    In Section \ref{sec:lmi_ansatz} we briefly summarize the approach presented in \cite{Herty2023a}, followed by Section \ref{sec:ssc_ansatz} where the key points of \cite{Yang2024} are given.
    Then in Section \ref{sec:comparison} we state the main result of this work which identifies SSC systems as particular representatives that satisfy an associated LMI.
    The obtained results are emphasized by discussing the example for the Saint-Venant equations, also presented in \cite{Yang2024}, in Section \ref{sec:ex_saintvenant}. To complete the discussion, we present in Section \ref{sec:diag_sys} an example where both approaches given in \cite{Yang2024,Herty2023} fail, but the approach presented in \cite{Herty2023} can be applied. Although \cite{Herty2023} is not within the main focus of the present work we want to highlight that it covers an important class of systems not covered by \cite{Herty2023a, Yang2024}.
    Numerical simulations for both examples are given in Section \ref{sec:num_res}.
    \section{The generic LMI approach}\label{sec:lmi_ansatz}
    In \cite{Herty2023a} the subsequent system of hyperbolic PDEs is studied
    \begin{align}
        \p{}{t}\mathbf{w}(t,\mathbf{x}) + \sum_{k=1}^d\mathbf{A}^{(k)}(\mathbf{x})\p{}{x_k}\mathbf{w}(t,\mathbf{x}) + \mathbf{B}(\mathbf{x})\mathbf{w}(t,\mathbf{x}) &= 0,\;(t,\mathbf{x})\in[0,T)\times\Omega
        \label{eq:hyp_cons_sys1}
    \end{align}
    Here $\mathbf{w}(t,\mathbf{x}) \equiv (w_1(t,\mathbf{x}),\dots,w_n(t,\mathbf{x}))^T$ is the vector of unknowns
    and $\Omega \subset \R^{d}$ a bounded domain with piecewise $C^1$ smooth boundary $\partial\Omega$.
    Moreover, $\mathbf{A}^{(k)}(\mathbf{x})$ and $\mathbf{B}(\mathbf{x})$ are sufficiently smooth and bounded $n\times n$ real matrices. The $\mathbf{A}^{(k)}(\mathbf{x})$ are in particular assumed to be symmetric.
    Usually symmetric hyperbolic systems appear with a general strictly positive definite symmetric matrix $\mathbf{A}^{(0)}$ in front of the time derivative.
    However, upon applying a suited variable transformation we may transform a symmetric hyperbolic system into the given form with the identity in front of the time derivative, cf. \cite{BenzoniGavage2007}.
    The assumption of symmetry is no major restriction since it includes all systems equipped with an additional conservation law, cf. \cite{Friedrichs1971,Dafermos2016}.
    This includes most systems relevant for applications, see \cite{Boillat1996,Ruggeri2021}.\\
    It is further assumed that there exists a \emph{feasible Lyapunov potential} $\overline{\mu}(\mathbf{x})$ such that
    \begin{align}
        \overline{\mathbf{m}} := \nabla\overline{\mu}(\mathbf{x})\quad\text{and}\quad \overline{\mathbf{A}}(\overline{\mathbf{m}}) := -\mathbf{Id} + \sum_{k=1}^d\overline{m}_k\mathbf{A}^{(k)} \geq 0.
        \label{cond:feas_potential}
    \end{align}
    Then
    \begin{align}
        \mathbf{A}(\mathbf{m}) := C\mathbf{Id} + \sum_{k=1}^dm_k\mathbf{A}^{(k)} \leq 0,\, C \in \R_{>0}\label{ineq:control_LMI}
    \end{align}
    holds for system \eqref{eq:hyp_cons_sys1} with $\mu(\mathbf{x}) = -C\overline{\mu}(\mathbf{x})$ and $\mathbf{m} = \nabla\mu(\mathbf{x})$.
    It is remarked in \cite{Herty2023a} that the LMI \eqref{ineq:control_LMI} can be modified if certain reminder terms, such as the coupling matrix $\mathbf{B}$, should be taken into account.
    Therefore we introduce with $\mathbf{B}^{sym} = \frac{1}{2}\left(\mathbf{B} + \mathbf{B}^T\right)$
    \[
      \mathcal{R}(\mathbf{x}) := \sum_{k=1}^d\p{}{x_k}\mathbf{A}^{(k)}(\mathbf{x}) - 2\mathbf{B}^{sym}(\mathbf{x})
    \]
    and demand
    \begin{align}
        \mathbf{A}(\mathbf{m}) := C\mathbf{Id} + \mathcal{R}(\mathbf{x}) + \sum_{k=1}^dm_k(\mathbf{x})\mathbf{A}^{(k)}(\mathbf{x}) \leq 0,\, C \in \R_{>0}.\label{ineq:control_LMI2}
    \end{align}
    The LMI \eqref{ineq:control_LMI2} then replaces \eqref{ineq:control_LMI}.
    For example in the case of  $\mathcal{R}(\mathbf{x}) < 0$ we could benefit from this additional term in order to find suited coefficients $\mathbf{m}$ as will be demonstrated below.
    The Lyapunov function is defined as follows
    \begin{align}
        L_1(t) = \int_\Omega \mathbf{w}(t,\mathbf{x})^T\mathbf{w}(t,\mathbf{x})\exp(\mu(\mathbf{x}))\,\dd\mathbf{x}.\label{eq:lyapunov_gen}
    \end{align}
    It is then shown, that under suited boundary conditions the Lyapunov function decays exponentially, see \cite{Herty2023a}[Thm. 2.4].
    \section{The approach for SSC systems}\label{sec:ssc_ansatz}
    In the following we briefly recall the main assumptions given in \cite{Yang2024} and we refer to this work for further details.
    In \cite{Yang2024} linear hyperbolic systems with constant coefficients in multi -- dimensions which satisfy the SSC are considered.
    We refer to the aforementioned references for further reading.
    It is shown that such a system can be written in the following form
    \begin{align}
        \p{}{t}\mathbf{U}(t,\mathbf{x}) + \sum_{k=1}^d\overline{\mathbf{A}}^{(k)}\p{}{x_k}\mathbf{U}(t,\mathbf{x}) = \overline{\mathbf{B}}\mathbf{U}(t,\mathbf{x}),\;(t,\mathbf{x})\in[0,T)\times\Omega
        \label{eq:hyp_cons_sys2}
    \end{align}
    with $\mathbf{U} = (\overline{\mathbf{u}},\overline{\mathbf{q}})^T$ where $\overline{\mathbf{u}} \in \R^{n-r}$ and $\overline{\mathbf{q}} \in \R^r$, with $0 < r < n$.
    Furthermore we have
    \[
      \overline{\mathbf{B}} = \begin{pmatrix}\mathbf{0}_{(n-r)\times (n-r)} & \mathbf{0}_{(n-r)\times r}\\ \mathbf{0}_{r\times (n-r)} & \overline{\mathbf{e}}\end{pmatrix}
    \]
    with $\overline{\mathbf{e}} \in \R^{r\times r}$ being invertible. The Jacobians are given by
    \[
      \overline{\mathbf{A}}^{(k)} = \begin{pmatrix}\overline{\mathbf{a}}_k & \overline{\mathbf{b}}_k\\ \overline{\mathbf{c}}_k & \overline{\mathbf{d}}_k\end{pmatrix},\,k=1,\dots,d
    \]
    with $\overline{\mathbf{a}}_k,\overline{\mathbf{b}}_k\in\R^{(n-r)\times(n-r)}$ and $\overline{\mathbf{c}}_k,\overline{\mathbf{d}}_k\in\R^{r\times r}$.
    In \cite{Yang2024}, system \eqref{eq:hyp_cons_sys2} is assumed to have the following properties
    \begin{enumerate}[(i)]
        \item There exists a symmetric positive definite matrix
        \[
          \mathbf{A}^{(0)} = \begin{pmatrix}\mathbf{X}_1 & \mathbf{0}_{(n-r)\times r}\\ \mathbf{0}_{r\times(n-r)} & \mathbf{X}_2\end{pmatrix}
        \]
        with $\mathbf{X}_1 \in \R^{(n-r)\times (n-r)}$ and $\mathbf{X}_2 \in \R^{r\times r}$, such that all matrices $\mathbf{A}^{(0)}\overline{\mathbf{A}}^{(k)}$ are symmetric.
        \item $\mathbf{X}_2\overline{\mathbf{e}} + \overline{\mathbf{e}}^T\mathbf{X}_2$ is positive definite.
        \item There exist real numbers $\alpha_k, k=1,\dots,d$ such that
        \[
          \sum_{k=1}^d \alpha_k\overline{\mathbf{a}}_k \in \R^{(n-r)\times(n-r)}
        \]
        has only negative eigenvalues.
    \end{enumerate}
    Following \cite{Yang2024} the properties (i) and (ii) are implied by the SSC.
    The Lyapunov function is then defined as follows
    \begin{align}
        L_2(t) = \int_\Omega \lambda(\mathbf{x})\mathbf{U}(t,\mathbf{x})^T\mathbf{A}^{(0)}\mathbf{U}(t,\mathbf{x})\,\dd\mathbf{x}.\label{eq:lyapunov_yy}
    \end{align}
    with $\lambda(\mathbf{x}) = K + \sum_{k=1}^d\alpha_k x_k > 0$.
    Note that the positivity of the weight function $\lambda(\mathbf{x})$ is an issue needed to be resolved by $K > 0$ with respect to domain $\Omega$.
    It is then shown in \cite{Yang2024}[Thm. 3.1, Lem. 3.2], that under suited boundary conditions the Lyapunov function decays exponentially.
    \section{Main result}\label{sec:comparison}
    In the subsequent part we state the main result of this work exploiting the relation between the approaches presented in \cite{Herty2023a,Yang2024}.
    \begin{thm}\label{thm:main}
    	Let system \eqref{eq:hyp_cons_sys2} be given and assume it satisfies the SSC property, i.e. the properties (i) - (iii) hold.
    	Then there exists a feasible Lyapunov potential $\mu(\mathbf{x})$ such that the LMI \eqref{ineq:control_LMI2} holds and the Lyapunov function \eqref{eq:lyapunov_gen} decays exponentially.
    \end{thm}
    \begin{prf}
    %
    %Property (i)
    %
    Property (i) states the existence of a symmetrizer for system \eqref{eq:hyp_cons_sys2}. In particular $\mathbf{X}_1$ and $\mathbf{X}_2$ are positive definite symmetric matrices.
    Multiplying \eqref{eq:hyp_cons_sys2} from the left with $\mathbf{A}^{(0)}$ gives 
    \begin{align*}
        \mathbf{A}^{(0)}\p{}{t}\mathbf{U}(t,\mathbf{x}) + \sum_{k=1}^d\mathbf{A}^{(0)}\overline{\mathbf{A}}^{(k)}\p{}{x_k}\mathbf{U}(t,\mathbf{x})
        &= \mathbf{A}^{(0)}\overline{\mathbf{B}}\mathbf{U}(t,\mathbf{x})\\
        \Leftrightarrow:\quad\mathbf{A}^{(0)}\p{}{t}\mathbf{U}(t,\mathbf{x}) + \sum_{k=1}^d\tilde{\mathbf{A}}^{(k)}\p{}{x_k}\mathbf{U}(t,\mathbf{x}) &= \tilde{\mathbf{B}}\mathbf{U}(t,\mathbf{x})
    \end{align*}
    with $\tilde{\mathbf{A}}^{(k)}, k=1,
    \dots,d$ being symmetric.
    Introducing the variables $\mathbf{w} := (\mathbf{A}^{(0)})^{\frac{1}{2}}\mathbf{U}$, $\mathbf{w} \equiv (\mathbf{u},\mathbf{q})^T$
    and multiplying from left by $(\mathbf{A}^{(0)})^{-\frac{1}{2}}$ we can transform the system to
    \begin{align*}
        &\p{}{t}\mathbf{w}(t,\mathbf{x}) + \sum_{k=1}^d\mathbf{A}^{(k)}\p{}{x_k}\mathbf{w}(t,\mathbf{x}) = -\mathbf{B}\mathbf{w}(t,\mathbf{x})\\
        \text{with}\quad &\mathbf{A}^{(k)} := (\mathbf{A}^{(0)})^{-\frac{1}{2}}\tilde{\mathbf{A}}^{(k)}(\mathbf{A}^{(0)})^{-\frac{1}{2}}\quad\text{being symmetric}\\
        \text{and}\quad &\mathbf{B} := -(\mathbf{A}^{(0)})^{-\frac{1}{2}}\tilde{\mathbf{B}}(\mathbf{A}^{(0)})^{-\frac{1}{2}}
    \end{align*}
    The transforms are directly applied to the respective block matrices, i.e.
    \begin{align*}
        \mathbf{a}_k &:= \mathbf{X}_1^{\frac{1}{2}}\overline{\mathbf{a}}_k\mathbf{X}_1^{-\frac{1}{2}},\quad
        \mathbf{b}_k  := \mathbf{X}_1^{\frac{1}{2}}\overline{\mathbf{b}}_k\mathbf{X}_1^{-\frac{1}{2}},\quad
        \mathbf{c}_k  := \mathbf{X}_2^{\frac{1}{2}}\overline{\mathbf{c}}_k\mathbf{X}_2^{-\frac{1}{2}},\quad
        \mathbf{d}_k  := \mathbf{X}_2^{\frac{1}{2}}\overline{\mathbf{d}}_k\mathbf{X}_2^{-\frac{1}{2}}\\
        \text{and}\quad \mathbf{e} &:= \mathbf{X}_2^{\frac{1}{2}}\overline{\mathbf{e}}\mathbf{X}_2^{-\frac{1}{2}}.
    \end{align*}
    The stated properties (ii) and (iii) are also transferred to the transformed system, which can be seen as follows.
    Property (ii) states that $\mathbf{X}_2\overline{\mathbf{e}} + \overline{\mathbf{e}}^T\mathbf{X}_2 > 0$. From this we yield for $\mathbf{v}\in\R^r\setminus\{\mathbf{0}\}$
    \begin{align*}
        0 &< \mathbf{v}^T\left(\mathbf{X}_2\overline{\mathbf{e}} + \overline{\mathbf{e}}^T\mathbf{X}_2\right)\mathbf{v}\\
        &= \mathbf{v}^T\left(\mathbf{X}_2\overline{\mathbf{e}}\mathbf{X}^{-\frac{1}{2}}_2\mathbf{X}^{\frac{1}{2}}_2
        + \mathbf{X}^{\frac{1}{2}}_2\mathbf{X}^{-\frac{1}{2}}_2\overline{\mathbf{e}}^T\mathbf{X}_2\right)\mathbf{v}\\
        &= \mathbf{v}^T\mathbf{X}_2^{\frac{1}{2}}\left(\mathbf{X}^{\frac{1}{2}}_2\overline{\mathbf{e}}\mathbf{X}^{-\frac{1}{2}}_2 +
        \mathbf{X}^{-\frac{1}{2}}_2\overline{\mathbf{e}}^T\mathbf{X}^{\frac{1}{2}}_2\right)\mathbf{X}^{\frac{1}{2}}_2\mathbf{v}\\
        &= \mathbf{y}^T\left(\mathbf{e} + \mathbf{e}^T\right)\mathbf{y}\quad\text{with}\quad\mathbf{y} = \mathbf{X}^{\frac{1}{2}}_2\mathbf{v}.
    \end{align*}
    Thus for the transformed system property (ii) states $\mathbf{e} + \mathbf{e}^T > 0$.\\
    According to property (iii) there exist real numbers $\alpha_k$ such that
    \[
      \overline{\mathbf{a}} := \sum_{k=1}^d\alpha_k\overline{\mathbf{a}}_k \in \R^{(n-r)\times(n-r)}
    \]
    has only negative eigenvalues. We now study the similarity transform
    \begin{align*}
        \mathbf{X}_1^{\frac{1}{2}}\overline{\mathbf{a}}\mathbf{X}_1^{-\frac{1}{2}} &= \mathbf{X}_1^{\frac{1}{2}}\left(\sum_{k=1}^d\alpha_k\overline{\mathbf{a}}_k\right)\mathbf{X}_1^{-\frac{1}{2}}
        = \sum_{k=1}^d\alpha_k\mathbf{X}_1^{\frac{1}{2}}\overline{\mathbf{a}}_k\mathbf{X}_1^{-\frac{1}{2}}
        = \sum_{k=1}^d\alpha_k\mathbf{a}_k =: \mathbf{a}.
    \end{align*}
    Since a similarity transform leaves the eigenvalues unchanged, $\mathbf{a}$ has only negative real eigenvalues.
    In summary we can state, that property (i) allows us to write the studied system \eqref{eq:hyp_cons_sys2} in the form
    \begin{align}
        \p{}{t}\mathbf{w}(t,\mathbf{x}) + \sum_{k=1}^d\mathbf{A}^{(k)}\p{}{x_k}\mathbf{w}(t,\mathbf{x}) + \mathbf{B}\mathbf{w}(t,\mathbf{x}) = 0.\label{symm_hyp_sys}
    \end{align}
    This is in particular a symmetric hyperbolic system as given by \eqref{eq:hyp_cons_sys1}.\\
    %
    %Property (ii) & (iii)
    %
    We study property (ii) and (iii) in terms of the transformed symmetric system \eqref{symm_hyp_sys}, i.e., $\mathbf{e} + \mathbf{e}^T > 0$ and
    $\mathbf{a}$ has only negative real eigenvalues for suited $\alpha_k\in\R$.
    To understand the implication of property (ii) and (iii) for the LMI 
    we differentiate the Lyapunov function with respect to time and rearrange the terms,
    such that the derivative can be written as the sum of a boundary integral $\mathcal{B}(t)$ and a volume integral $\mathcal{I}(t)$.
    Since \eqref{eq:hyp_cons_sys2} has constant coefficient matrices we obtain
    \begin{align*}
        \frac{\dd}{\dd t}L_1(t) &= -\underbrace{\int_{\partial\Omega} \mathcal{A}(t,\mathbf{x})\cdot\mathbf{n}(\mathbf{x})\,\dd s(\mathbf{x})}_{=:\mathcal{B}(t)}\\
        &+ \underbrace{\int_\Omega \mathbf{w}^T(t,\mathbf{x})\left[\sum_{k=1}^dm_k\mathbf{A}^{(k)} - 2\mathbf{B}\right]\mathbf{w}(t,\mathbf{x})\exp(\mu(\mathbf{x}))\,\dd\mathbf{x}}_{=: \mathcal{I}(t)}.
    \end{align*}
    Without loss of generality we replace $\mathbf{B}$ by $\mathbf{B}^{sym} = \left(\mathbf{B} + \mathbf{B}^T\right)/2$.
    We now have to show that the LMI \eqref{ineq:control_LMI2} holds. Thus we use $\mathbf{w} = (\mathbf{u}, \mathbf{q})^T$ and yield the following quadratic form
    \begin{align*}
        \mathbf{w}^T\mathbf{A}^{(k)}\mathbf{w} = \mathbf{u}^T\mathbf{a}_k\mathbf{u} + \mathbf{u}^T\mathbf{b}_k\mathbf{q}
        + \mathbf{q}^T\mathbf{c}_k\mathbf{u} + \mathbf{q}^T\mathbf{d}_k\mathbf{q}.
    \end{align*}
    Further the term
    \[
      -2\mathbf{w}^T\mathbf{B}\mathbf{w} =  -2\mathbf{w}^T\mathbf{B}^{sym}\mathbf{w} = -2\mathbf{q}^T\left(\mathbf{e} + \mathbf{e}^T\right)\mathbf{q}
    \]
    is scaled with a constant $K > 0$.
    This gives
    \[
      \mathbf{A}(\mathbf{m}) := C\mathbf{Id} - 2\mathbf{B}^{sym} + \sum_{k=1}^dm_k\mathbf{A}^{(k)} = C\mathbf{Id} + \frac{1}{K}\left(-2K\mathbf{B}^{sym} + \sum_{k=1}^d\tilde{m}_k\mathbf{A}^{(k)}\right).
    \]
    with $\tilde{m}_k = Km_k$. Due to properties (ii) and (iii) by setting $\tilde{m}_k = \alpha_k$ we then yield for system \eqref{symm_hyp_sys}
    \begin{align*}
        \mathbf{w}^T\mathbf{A}(\mathbf{m})\mathbf{w} &= C\|\mathbf{w}\|_2^2 + \frac{1}{K}\left(-K\mathbf{q}^T\left(\mathbf{e} + \mathbf{e}^T\right)\mathbf{q}^T
        + \sum_{k=1}^d\alpha_k\mathbf{w}^T\mathbf{A}^{(k)}\mathbf{w}^T\right)\\
        &= C\|\mathbf{w}\|_2^2 + \frac{1}{K}\big(-K\mathbf{q}^T\left(\mathbf{e} + \mathbf{e}^T\right)\mathbf{q}^T
        + \underbrace{\mathbf{u}^T\left(\sum_{k=1}^d\alpha_k\mathbf{a}_k\right)\mathbf{u}}_{\stackrel{(ii)}{<}0}\\
        &\left.+ \mathbf{u}^T\left(\sum_{k=1}^d\alpha_k\mathbf{b}_k\right)\mathbf{q}
        + \mathbf{q}^T\left(\sum_{k=1}^d\alpha_k\mathbf{c}_k\right)\mathbf{u} + \mathbf{q}^T\left(\sum_{k=1}^d\alpha_k\mathbf{d}_k\right)\mathbf{q}\right)\\
    \end{align*}
    Now due to $\mathbf{e} + \mathbf{e}^T > 0$ the constant $K > 0$ is chosen such that
    \[
      \mathbf{u}^T\left(\sum_{k=1}^d\alpha_k\mathbf{b}_k\right)\mathbf{q}
        + \mathbf{q}^T\left(\sum_{k=1}^d\alpha_k\mathbf{c}_k\right)\mathbf{u} + \mathbf{q}^T\left(\sum_{k=1}^d\alpha_k\mathbf{d}_k\right)\mathbf{q} -K\mathbf{q}^T\left(\mathbf{e} + \mathbf{e}^T\right)\mathbf{q}^T < 0
    \]
    Thus there exists a constant $C > 0$ such that \eqref{ineq:control_LMI2} holds with $m_k = \alpha_k/K$. Applying \cite{Herty2023a}[Thm. 2.4] finishes the proof.
    \end{prf}
    \begin{rem}
    	Within the realm of Thm. \ref{thm:main} both approaches give exponential decay regardless of the weight used for the Lyapunov function. This non-uniqueness is also mentioned in \cite{Yang2024}.
    	Nevertheless, for the sake of completeness we want to comment on the Lyapunov function in more detail and in particular the different weight functions.
	    The Lyapunov function used in the framework of \cite{Yang2024} is given as $\lambda(x,y) = K + \alpha x + \beta y$. Following \cite{Herty2023a} we have for the situation under study
	    \[
	      f(x,y) := \exp(\mu(x,y))\quad\text{with}\quad\mu(x,y) = m_1x + m_2y = \frac{\alpha}{K}x + \frac{\beta}{K}y + \ln(K)
	    \]
	    The additional term $\ln(K)$ is just added for proper scaling and does not affect the results of \cite{Herty2023a}.
	    For the exponential weight function we thus obtain by a convexity argument 
	    \begin{align*}
	        f(x,y) = \exp\left(\frac{\alpha}{K}x + \frac{\beta}{K}y + \ln(K)\right)
	        %&= f(x_0,y_0) + \nabla f(x_0,y_0)^T\begin{pmatrix}x-x_0\\y-y_0\end{pmatrix}% + \frac{1}{2}(x-x_0,y-y_0)\mathbf{D}^2f(\xi_x,\xi_y)\begin{pmatrix}x-x_0\\y-y_0\end{pmatrix}\\
	        %= \exp(\mu(0,0)) + (m_1, m_2)\exp(\mu(0,0))\begin{pmatrix}x\\y\end{pmatrix}\\
	        %&+ \underbrace{\frac{1}{2}(x,y)\begin{pmatrix} m_1^2 & m_1 m_2\\ m_1 m_2 &
	        %m_2^2\end{pmatrix}\exp(\mu(\xi_x,\xi_y))\begin{pmatrix}x\\y\end{pmatrix}}_{\geq 0}\\
	        \geq K + K\left(\frac{\alpha}{K}x + \frac{\beta}{K}y\right) = \lambda(x,y).
	    \end{align*}
	    Moreover, by choosing the weight functions as stated above the Lyapunov functions are related by
	    \begin{align}
	        \int_\Omega \mathbf{w}(t,x,y)^T\mathbf{w}(t,x,y)\exp(\mu(x,y))\,\dd\mathbf{x} \geq \int_\Omega \mathbf{w}(t,x,y)^T\mathbf{w}(t,x,y)\lambda(x,y)\,\dd\mathbf{x}.\label{ineq:lyapunov_comp}
	    \end{align}
	    This implies in the case of exponential decay in \cite{Herty2023a}
	    \begin{align*}
	        \frac{\dd}{\dd t}\int_\Omega \mathbf{w}(t,x,y)^T\mathbf{w}(t,x,y)\exp(\mu(x,y))\,\dd\mathbf{x}  &\leq -C\int_\Omega \mathbf{w}(t,x,y)^T\mathbf{w}(t,x,y)\exp(\mu(x,y))\,\dd\mathbf{x}\\ &\leq -C\int_\Omega \mathbf{w}(t,x,y)^T\mathbf{w}(t,x,y)\lambda(x,y)\,\dd\mathbf{x}.
	    \end{align*}
	\end{rem}
	\begin{rem}
	    Concerning the boundary condition there is no significant structural difference between the approaches presented in \cite{Yang2024} and \cite{Herty2023a}.
	    Due to the symmetry of the Jacobians $\mathbf{A}^{(k)}$ the pencil matrix of the system
	    \[
	      \mathbf{A}^\ast(\mathbf{\nu}) := \sum_{k=1}^d \nu_k\mathbf{A}^{(k)}
	    \]
	    is symmetric and therefore diagonalizable on the boundary. In particular we have
	    \[
	      \mathbf{w}^T\mathbf{A}^\ast(\mathbf{n})\mathbf{w} = \mathbf{w}^T\left(\sum_{k=1}^d n_k\mathbf{A}^{(k)}\right)\mathbf{w} = \mathbf{w}^T\mathbf{T}\mathbf{T}^T\left(\sum_{k=1}^d n_k\mathbf{A}^{(k)}\right)\mathbf{T}\mathbf{T}^T\mathbf{w} = \mathbf{v}^T\mathbf{\Lambda}\mathbf{v}
	    \]
	    where $\mathbf{n}$ denotes the outward pointing normal of $\partial\Omega$, $\mathbf{T}$ is the orthogonal transformation matrix,
	    $\mathbf{v} = \mathbf{T}^T\mathbf{w}$ and $\mathbf{\Lambda}$ is the diagonal matrix with the eigenvalues the system.
    \end{rem}
    \begin{rem}
    	The case of an arbitrary matrix $\mathbf{B}^{sym}$ is excluded in \cite{Yang2024} due to the SSC. For an example where the SSC is violated due to a regular right hand side we refer to Ssction \ref{sec:diag_sys}. In the case of an arbitrary right hand side we follow \cite{Herty2023,Herty2023a} and estimate the quadratic form as follows
    \[
      -2\mathbf{w}^T(t,x,y)\mathbf{B}\mathbf{w}(t,x,y) =  -2\mathbf{w}^T(t,x,y)\mathbf{B}^{sym}\mathbf{w}(t,x,y) \leq C_B\|\mathbf{w}(t,\mathbf{x})\|_2^2,\,C_B\in\R_{>0}.
    \]
    This term then needs to be consumed by the remaining terms.

    Considering the case $\mathbf{B}^{sym} > 0$ we can simply estimate
    \[
      -2\mathbf{w}^T(t,x,y)\mathbf{B}\mathbf{w}(t,x,y) = -2\mathbf{w}^T(t,x,y)\mathbf{B}^{sym}\mathbf{w}(t,x,y) \leq 0.
    \]

    However, as noted before, it is also possible to benefit from the coupling term and relax \eqref{ineq:control_LMI} to \eqref{ineq:control_LMI2}.
    This is for example used in \cite{Yang2024} where the SSC holds and the idea is outlined in the presented proof of Thm. \ref{thm:main}.
    \end{rem}
    \section{Examples}
    \subsection{Application to the Saint-Venant system}\label{sec:ex_saintvenant}
    In \cite{Yang2024} an example for the Saint-Venant equations is presented and we want to give a short review here, showing that the approach given in \cite{Herty2023a} is indeed applicable.
    Note that in \cite{Herty2023a} an example for the barotropic Euler equations is given where the SSC does not hold.
    The linearized system in terms of the unknowns $\mathbf{w} = (\tilde{h}, w, v)^T$ is given by
    \begin{align}
        &\p{}{t}\mathbf{w}(t,x,y) + \mathbf{A}^{(1)}\p{}{x}\mathbf{w}(t,x,y) + \mathbf{A}^{(2)}\p{}{y}\mathbf{w}(t,x,y) = -\mathbf{B}\mathbf{w}(t,x,y)\label{eq:saintv_sys}\\
        \text{with}\quad &\mathbf{A}^{(1)} = \begin{pmatrix} w^\ast & \sqrt{gH^\ast} & 0\\ \sqrt{gH^\ast} & w^\ast & 0\\ 0 & 0 & w^\ast\end{pmatrix},\;
        \mathbf{A}^{(2)} = \begin{pmatrix} v^\ast & 0 & \sqrt{gH^\ast}\\ 0 & v^\ast & 0\\ \sqrt{gH^\ast} & 0 & v^\ast\end{pmatrix},\;
        \mathbf{B} = \begin{pmatrix} 0 & 0 & 0\\ 0 & k & -l\\ 0 & l & k\end{pmatrix}.\notag
    \end{align}
    where $\tilde{h} = \sqrt{\frac{g}{H^\ast}}h$ is a scaled perturbed height, $w$ the perturbed velocity in $x$-direction and $v$ the perturbed velocity in $y$-direction.
    The quantity $g$ is the gravitational acceleration, $k>0$ is the viscous drag coefficient and $l > 0$ is the Coriolis coefficient.
    Further $(H^\ast,w^\ast,v^\ast)^T$ is a steady state in terms of the height, the velocity in $x$-direction and the velocity in $y$-direction, respectively.
    For $\mathbf{B}$ we have that $\mathbf{B}^{sym} > 0$ which in the approach given in \cite{Herty2023a} could be estimated by zero.
    However, as remarked before we can also make use of $\mathbf{B}^{sym}$ for the estimate as in \eqref{ineq:control_LMI2} with $\tilde{C} > 0$.
    We reformulate the inequality to introduce a further scaling for $\mathbf{B}^{sym}$ with $\chi > 0$
    \begin{align*}
        0\geq \tilde{\mathbf{A}}(\mathbf{m}) &:= \tilde{C}\mathbf{Id} - 2\mathbf{B}^{sym} + \sum_{k=1}^d\tilde{m}_k\mathbf{A}^{(k)}\\
        &= \frac{1}{\chi}\left(\chi \tilde{C}\mathbf{Id} - 2\chi\mathbf{B}^{sym} + \sum_{k=1}^d\chi\tilde{m}_k\mathbf{A}^{(k)}\right)\\
        \Leftrightarrow\quad
        0 \geq \mathbf{A}(\mathbf{m}) &:= C\mathbf{Id} - 2\chi\mathbf{B}^{sym} + \sum_{k=1}^dm_k\mathbf{A}^{(k)}
    \end{align*}
    We have
    \[
      \mathbf{A}(\mathbf{m}) = \begin{pmatrix} C + m_1w^\ast + m_2v^\ast & m_1\sqrt{gH^\ast} & m_2\sqrt{gH^\ast}\\
      m_1\sqrt{gH^\ast} & C - 4\chi k + m_1w^\ast + m_2v^\ast & 0\\
      m_2\sqrt{gH^\ast} & 0 & C - 4\chi k + m_1w^\ast + m_2v^\ast \end{pmatrix}.
    \]
    In \cite{Yang2024} the situation is considered with the vertical velocity in the steady state $v^\ast=0$
    and with the horizontal velocity in the steady state $0 < w^\ast < \sqrt{gH^\ast}$, respectively. Thus we have with $m := m_1$
    \[
      \mathbf{A}(m) = \begin{pmatrix} C + m w^\ast & m\sqrt{gH^\ast} & 0\\
      m\sqrt{gH^\ast} & C - 4\chi k + m w^\ast & 0\\
      0 & 0 & C - 4\chi k + m w^\ast \end{pmatrix}.
    \]
    We obtain the following conditions on $m$ and $\chi$ to have $\mathbf{A}(m) \leq 0$
    \begin{enumerate}[(i)]
        \item $C + m w^\ast < 0\quad\Leftrightarrow\quad m < -\frac{C}{w^\ast}$
        \item Due to $4\chi k > 0$ condition (i) also implies $0 > C + m w^\ast > C - 4\chi k +  mw^\ast$
        \item For the second principle minor we yield
        \[
          \left(C + m w^\ast\right)\left(C - 4\chi k + mw^\ast\right) - m^2gH^\ast = \left(C + m w^\ast\right)^2 - 4\chi k\underbrace{\left(C + m w^\ast\right)}_{\leq 0} - m^2 gH^\ast
        \]
        and thus there exists a $\chi > 0$ such that this expression becomes positive, i.e.,
        \[
          \chi \geq \frac{\left(C + m w^\ast\right)^2 - m^2 gH^\ast}{4 k\left(C + m w^\ast\right)}.
        \]
        Hence for $m < -\frac{C}{w^\ast}$ and $\chi > 0$ large enough $\mathbf{A}(m) \leq 0$ holds.
        In the case of $m = -1$ and $\chi = 2L$, i.e., the parameters used in \cite{Yang2024}, these conditions impose a restriction on the decay rate, i.e.
        \begin{align}
        	\underbrace{w^\ast + 8Lk}_{(ii)} > \underbrace{w^\ast}_{(i)} > \underbrace{w^\ast + 4Lk - \sqrt{16L^2k^2 + gH^\ast}}_{(iii)} \geq C > 0.\label{constr_decay_rate}
        \end{align}
    \end{enumerate}
    It remains to study the boundary term and prescribe boundary conditions, such that $\mathcal{BC} \geq 0$.
    To this end we will make use of the eigenstructure of the system given in the Appendix \ref{app:eigenstruct} for the readers convenience.
    Following the example of \cite{Yang2024} the boundary of the domain $\Omega$ is given by $\partial\Omega = [0,L]\times\{0\}\cup\{L\}\times[0,1]\cup[0,L]\times\{1\}\cup\{0\}\times[0,1]$.
    For the boundary integral we have
    \begin{align*}
        \mathcal{BC} &= \int_{\partial\Omega}\mathbf{w}^T\mathbf{A}^\ast(\mathbf{n})\mathbf{w}\delta(x,y)\,\dd s(x,y)
        = \int_{\partial\Omega}\mathbf{v}^T\mathbf{\Lambda}^\ast(\mathbf{n})\mathbf{v}\delta(x,y)\,\dd s(x,y)\\
        &= \sum_{i=1}^3\int_{\partial\Omega}v_i^2\lambda_i(\mathbf{n})\delta(x,y)\,\dd s(x,y)
    \end{align*}
    where $\mathbf{v} = \mathbf{T}^T(\mathbf{n})\mathbf{w}$ is calculated according to \eqref{sv_vec_trafo} and $\mathbf{n}$ denotes the outward pointing normal of the boundary $\partial\Omega$.
    The weight function is given as follows
    \begin{align}
    	\delta(x,y) = \begin{dcases}
    		&\exp\left(\ln(2L) - \dfrac{x}{2L}\right),\,\text{\cite{Herty2023a}}\\
    		&2L - x, \, \text{\cite{Yang2024}}
    	\end{dcases}.\label{def:weight_fun}
    \end{align}
    Now we identify the controllable and uncontrollable parts of the boundary, i.e. for $i = 1,2,3$
    \begin{align*}
        \Gamma_i^+ &:= \left\{\left.\mathbf{x} \in \partial\Omega\,\right|\,\lambda_i(\mathbf{x},\mathbf{n}(\mathbf{x})) \geq 0\right\},\\
        \Gamma_i^- &:= \left\{\left.\mathbf{x} \in \partial\Omega\,\right|\,\lambda_i(\mathbf{x},\mathbf{n}(\mathbf{x})) < 0\right\}.
    \end{align*}
    These are given by
    \begin{align}
        \Gamma_1^- &= \partial\Omega,\quad &\Gamma_1^+ &= \emptyset,\notag\\
        \Gamma_2^- &= \{0\}\times(0,1),\quad &\Gamma_2^+ &= \partial\Omega\setminus\Gamma_2^-,\label{bound_split}\\
        \Gamma_3^- &= \emptyset,\quad &\Gamma_3^+ &= \partial\Omega.\notag
    \end{align}
    Let us denote the general controls by $\varphi_1(t,x,y)$ for the first component and $\psi_2(t,y)$ for the second component, respectively.
    These have to be chosen such that
    \begin{align*}
        \mathcal{BC} &= \underbrace{\int_{\partial\Omega}\varphi_1(t,x,y)^2(\mathbf{n}(x,y)w^\ast - \sqrt{gH^\ast})\delta(x,y)\,\dd s(x,y) - w^\ast\int_0^1\psi_2(t,y)^2\delta(0,y)\,\dd y}_{\leq 0}\\
        &+ \underbrace{w^\ast\int_{\partial\Omega\setminus\Gamma_2^-}v^2_2\mathbf{n}(x,y)\delta(x,y)\,\dd s(x,y)
        + \int_{\partial\Omega}v^2_3(\mathbf{n}(x,y)w^\ast + \sqrt{gH^\ast})\delta(x,y)\,\dd s(x,y)}_{\geq 0} \geq 0.
    \end{align*}
    A detailed derivation of the boundary conditions is given in the Appendix \ref{app:bound_contr} and numerical simulations are given below.
    Note that the homogeneous part of system \eqref{eq:saintv_sys} is exactly the system studied in \cite{Dia2013}. In \cite{Dia2013} a Lyapunov function of the form
    \begin{align}
        L(t) = \int_\Omega\mathbf{w}^T\mathcal{L}\mathbf{w}\,\dd\mathbf{x}\label{eq:lyapunov_dia}
    \end{align}
    is used with $\mathcal{L} := \diag(1,H^\ast/g,H^\ast/g)$. It is shown under suited assumptions that the system can be stabilized near $\mathbf{0}$ in terms of $L(t) \leq CL(0)$ with $C > 0$.
    So while in \cite{Dia2013} the boundedness of the Lyapunov function is shown, the approaches presented in \cite{Herty2023a} and \cite{Yang2024} are able to establish exponential decay.
    It should be remarked that in \cite{Dia2018} the same authors established exponential decay in a particular two dimensional situation and therefore used a reduced one dimensional system.
    Further note the approaches presented in \cite{Herty2023a} and \cite{Yang2024} are also applicable to more complex geometries.
    The choice of the rectangular domain here is due to \cite{Yang2024} and just for instructive purposes.
    \subsection{The case of a diagonal system}\label{sec:diag_sys}
    With the subsequent example we want to show where the approaches of \cite{Herty2023a} and \cite{Yang2024} are not applicable, but which can be treated by \cite{Herty2023}.
    Thus we want to emphasize that the case of systems diagonal Jacobians is not just a special case of symmetric systems, but a relevant class on its own.
    Therefore we now consider the following inhomogeneous system of three equations in two dimensions with constant Jacobians
    \begin{align}
    	&\p{}{t}\mathbf{w}(t,x,y) + \mathbf{A}^{(1)}\p{}{x}\mathbf{w}(t,x,y) + \mathbf{A}^{(2)}\p{}{y}\mathbf{w}(t,x,y) = -\mathbf{B}\mathbf{w}(t,x,y)\label{eq:diag_sys_ex}\\
        \text{with}\quad &\mathbf{A}^{(1)} = \begin{pmatrix} 1 & 0 & \hphantom{-}0\\ 0 & 1 & \hphantom{-}0\\ 0 & 0 & -1\end{pmatrix},\;
        \mathbf{A}^{(2)} = \begin{pmatrix} 1 & \hphantom{-}0 & \hphantom{-}0\\ 0 & -1 & \hphantom{-}0\\ 0 & \hphantom{-}0 & -1\end{pmatrix},\;
        \mathbf{B} = \begin{pmatrix} -1 & \hphantom{-}0 & \hphantom{-}0\\ 0 & \hphantom{-}2 & -1\\ 0 & -1 & \hphantom{-}1 \end{pmatrix}.\notag
    \end{align}
    Clearly the LMI approach for the homogeneous part can not be applied since we have three inequalities for two unknowns and their solution set is empty. Moreover, the matrix $\mathbf{B}$ is symmetric, indefinite and therefore it cannot be used in order help satisfying the LMI \eqref{ineq:control_LMI2}.
    Furthermore the matrix $\mathbf{B}$ is regular and thus the system does not satisfy the SSC. This can easily be seen since as part of the SSC there should exist a regular matrix $\mathbf{P}$ such that
    \[
    	\mathbf{P}\mathbf{B}\mathbf{P}^{-1} = \begin{pmatrix}\mathbf{0}_{(n-r)\times (n-r)} & \mathbf{0}_{(n-r)\times r}\\ \mathbf{0}_{r\times (n-r)} & \overline{\mathbf{e}}\end{pmatrix}.
    \]
    Obviously the RHS is singular and thus by the determinant product rule and due to $\mathbf{B}$ regular there exists no such $\mathbf{P}$.
    Since both approaches are not possible we want to demonstrate how this system is stabilized using the results given in \cite{Herty2023}.
    There the Lyapunov function is defined by
    \begin{align}
        L(t) = \int_\Omega \mathbf{w}(t,\mathbf{x})^T\mathcal{E}(\mu(\mathbf{x}))\mathbf{w}(t,\mathbf{x})\,\dd\mathbf{x}\label{eq:lyapunov_gen_diag}
    \end{align}
    with
	\begin{align*}
		\mathcal{E}(\mu(\mathbf{x})) &:= \textup{diag}(\exp(\mu_1(\mathbf{x})),\dots,\exp(\mu_n(\mathbf{x}))).
	\end{align*}
    We consider the spatial domain $\Omega = [0,1]\times [0,1]$. From the matrices \eqref{eq:diag_sys_ex} we obtain the following vectors
	\begin{align}
	    \mathbf{a}_1 = (1,1),\quad\mathbf{a}_2 = (1,-1)\quad\text{and}\quad\mathbf{a}_3 = (-1,-1).\label{ray_vecs}
	\end{align}
	For the normal of the domain $\Omega$ we obtain
	\begin{align*}
	    \mathbf{n}(\mathbf{x}) =
	    \begin{cases}
	        (1,0),\,&\mathbf{x}\in\{1\}\times[0,1],\\
	        (0,1),\,&\mathbf{x}\in[0,1]\times\{1\},\\
	        (-1,0),\,&\mathbf{x}\in\{0\}\times[0,1],\\
	        (0,-1),\,&\mathbf{x}\in[0,1]\times\{0\}
	    \end{cases}
	\end{align*}
	and thus we yield for the associated products $\mathbf{a}_i\cdot\mathbf{n}$
	\begin{align}
	    \mathbf{a}_1\cdot\mathbf{n}(\mathbf{x}) &=
	    \begin{cases}
	        \hphantom{-}1,\,&\mathbf{x}\in\{1\}\times[0,1],\\
	        \hphantom{-}1,\,&\mathbf{x}\in[0,1]\times\{1\},\\
	        -1,\,&\mathbf{x}\in\{0\}\times[0,1],\\
	        -1,\,&\mathbf{x}\in[0,1]\times\{0\}
	    \end{cases},\quad
	    \mathbf{a}_2\cdot\mathbf{n}(\mathbf{x}) =
	    \begin{cases}
	        \hphantom{-}1,\,&\mathbf{x}\in\{1\}\times[0,1],\\
	        -1,\,&\mathbf{x}\in[0,1]\times\{1\},\\
	        -1,\,&\mathbf{x}\in\{0\}\times[0,1],\\
	        \hphantom{-}1,\,&\mathbf{x}\in[0,1]\times\{0\}
	    \end{cases}\notag\\
	    \text{and}\quad
	    \mathbf{a}_3\cdot\mathbf{n}(\mathbf{x}) &=
	    \begin{cases}
	        -1,\,&\mathbf{x}\in\{1\}\times[0,1],\\
	        -1,\,&\mathbf{x}\in[0,1]\times\{1\},\\
	        \hphantom{-}1,\,&\mathbf{x}\in\{0\}\times[0,1],\\
	        \hphantom{-}1,\,&\mathbf{x}\in[0,1]\times\{0\}
	    \end{cases}.\label{normal_products}
	\end{align}
	Hence we have the following partitioning of the boundary
	\begin{align}
	    \begin{split}\label{domain_split}
	        \Gamma_1^+ &:= \{1\}\times[0,1]\cup[0,1]\times\{1\},\quad\Gamma_1^- := \{0\}\times[0,1]\cup[0,1]\times\{0\},\\
	        \Gamma_2^+ &:= \{1\}\times[0,1]\cup[0,1]\times\{0\},\quad\Gamma_2^- := \{0\}\times[0,1]\cup[0,1]\times\{1\},\\
	        \Gamma_3^+ &:= \{0\}\times[0,1]\cup[0,1]\times\{0\},\quad\Gamma_3^- := \{1\}\times[0,1]\cup[0,1]\times\{1\}.
	    \end{split}
	\end{align}
	In what follows we apply the control for the first and second component on the part $\mathcal{C}_1 = \mathcal{C}_2 = \{0\}\times[0,1]$ and prescribe zero boundary conditions on $\mathcal{Z}_1 = [0,1]\times\{0\}$
	and $\mathcal{Z}_2 = [0,1]\times\{1\}$.
	For the third component we prescribe the control on $\mathcal{C}_3 = [0,1]\times\{1\}$ and prescribe zero boundary conditions on $\mathcal{Z}_3 = \{1\}\times[0,1]$. Note that the boundary parts correspond to the respective component of $\mathbf{w}$.
	Before we determine the weight functions $\mu_i$ we need to discuss the coupling matrix $\mathbf{B}$. Clearly $\mathbf{B}$ is symmetric and indefinite with the eigenvalues
	\[
	  \lambda_1 = -1 < 0 <\lambda_{2,3} = \frac{3\pm\sqrt{5}}{2}.
	\]
	Due to the symmetry we can establish the following estimate
	\begin{align*}
		-\mathbf{w}^T\left(\mathbf{B}^T\mathcal{E} + \mathcal{E}\mathbf{B}\right)\mathbf{w} \leq 2\mathbf{w}^T\mathcal{E}\mathbf{w}. 
	\end{align*}
	For further details and insight we refer to \cite{Herty2023}.
	Next we determine the weight functions $\mu_1(x,y),\mu_2(x,y)$ and $\mu_3(x,y)$. Due to the given structure we obtain for the weight functions
	\begin{align*}
	    &\mathbf{a}_i\cdot\nabla\mu_i(x,y) + 2 = -C_L^{(i)}
	    \quad\Leftrightarrow\quad \p{}{x}\mu_i(x,y) + \frac{\mathbf{a}^{(2)}_i}{\mathbf{a}_i^{(1)}}\p{}{y}\mu_i(x,y) = -\frac{C_L^{(i)} + 2}{\mathbf{a}_i^{(1)}}\\
	    \Rightarrow\quad &\mu_i(\mathbf{x}) = g_i\left(y - \frac{\mathbf{a}^{(2)}_i}{\mathbf{a}_i^{(1)}}x\right) -\frac{C_L^{(i)} + 2}{\mathbf{a}_i^{(1)}}x.
	\end{align*}
	This holds for arbitrary $g_i \in C^1$ and for the sake of simplicity we assume $g_i(\sigma) = \sigma$ and further $C_L^{(1)} = C_L^{(2)} = C_L^{(3)} = C_L > 0$.
	Thus we yield the following weight functions
	\begin{align}
	    \begin{split}\label{weight_fun}
	    	\mu_1(x,y) &= y - x - \left(C_L + 2\right)x,\quad\mu_2(x,y) = y + x - \left(C_L + 2\right)x\\
	    	\text{and}\quad\mu_3(x,y) &= y - x + \left(C_L + 2\right)x.
	    \end{split}
	\end{align}
	Next we specify the constraint for the boundary control and we assume $u(t) := u_1(t,x,y) = u_2(t,x,y) = u_3(t,x,y)$ which leads to
	\begin{align}
		u(t)^2 \leq -\left(\sum_{i=1}^3\int_{\mathcal{C}_i}\left(\mathbf{a}_i\cdot\mathbf{n}\right)\exp(\mu_i(x,y))\,\dd s(x,y)\right)^{-1}
		\sum_{i=1}^3\int_{\Gamma_i^+}w_i^2\left(\mathbf{a}_i\cdot\mathbf{n}\right)\exp(\mu_i(x,y))\,\dd s(x,y).\label{ineq:boundary_control2}
	\end{align}
	Inserting the obtained results we get
	\begin{align}
	    &-\sum_{i=1}^3\int_{\mathcal{C}_i}\left(\mathbf{a}_i\cdot\mathbf{n}\right)\exp(\mu_i(x,y))\,\dd s(x,y) = 2\int_0^1\exp(y)\,\dd y + \int_0^1\exp(1 + (C_L + 1)x)\,\dd x\notag\\
	    &= 2(e - 1) + \frac{e}{C_L + 1}(\exp(C_L + 1) - 1) =: C(C_L) > 0,\notag\\
	    \mathcal{I}(t) &:= \sum_{i=1}^3\int_{\Gamma_i^+}w_i^2\left(\mathbf{a}_i\cdot\mathbf{n}\right)\exp(\mu_i(\mathbf{x}))\,\dd\mathbf{x}\notag\\
	    &= \int_0^1 w_1(t,1,\sigma)^2\exp(\mu_1(1,\sigma)) + w_1(t,\sigma,1)^2\exp(\mu_1(\sigma,1))\dots\notag\\ &\dots + w_2(t,\sigma,0)^2\exp(\mu_2(\sigma,0)) + w_2(t,1,\sigma)^2\exp(\mu_2(1,\sigma))\dots\notag\\
	    &\dots + w_3(t,\sigma,0)^2\exp(\mu_3(\sigma,0)) + w_3(t,0,\sigma)^2\exp(\mu_3(0,\sigma))\,\dd\sigma\notag\\
	    \Rightarrow\quad &u(t)^2 \leq \frac{1}{C(C_L)}\mathcal{I}(t).\label{ineq:u_bc}
	\end{align}
	Thus a possible boundary control satisfying the premisses of the main theorem in \cite{Herty2023} is given by
	\begin{align}
		u(t) = \sqrt{\frac{1}{C(C_L)}\mathcal{I}(t)}.\label{eq:u_bc}
	\end{align}
	Now we have specified all the details needed to stabilize the solution of the given PDE \eqref{eq:diag_sys_ex} in $\Omega$. This will be demonstrated in the following section numerically.
	%
	%%%%%%%%%%%%%%%%%%%%%%%%%%%%%%
	%
	\section{Numerical Results}\label{sec:num_res}
	\subsection{The Saint-Venant system}\label{sec:sv_sys_num}
	In this section the control for the Saint-Venant is studied numerically. The MUSCL second--order finite--volume scheme, see \cite{Toro2009}, is used on a regular mesh for $\Omega = [0,L]\times[0,1]$ with grid size $\Delta x \times \Delta y$ to solve the discretized equation \eqref{eq:diag_sys_ex}.
	The cell average of $\mathbf{W}$ on $C_{i,j} = [x_i-\frac{\Delta x}2, x_i+\frac{\Delta x}2] \times [y_j-\frac{\Delta y}2, y_j+\frac{\Delta y}2]$, $x_i=i \Delta x, y_j = j\Delta y$ and time $t_n = n\Delta t$ is given by
	\begin{align*}
	    \mathbf{W}_{i,j}^n = \frac{1}{|C_{i,j}|} \int_{C_{i,j}} \mathbf{W}(t_n,x,y)\,\dd x\dd y,
	\end{align*}
	for $i=0,\dots,N_x, j=0,\dots,N_y$ and $n\geq 0$ and where $N_x \Delta x = L$ and $N_y\Delta y = 1$. The cell averages of the initial data $\mathbf{W}_0$ are obtained analogously and define $\mathbf{W}_{i,j}^0$.
	Transmissive boundary conditions are used for outgoing waves, see \cite{Toro2009}.
	The primitive boundary controls \eqref{contr_a_primvar_bp1} -- \eqref{contr_primvar_bp4} are imposed otherwise.
	The Lyapunov function \eqref{eq:lyapunov_gen} is approximated at time $t_n$ by $L^n$ using a numerical quadrature rule on the equi-distant grid 
	\begin{align}
	    L^n:=\Delta x\Delta y \sum\limits_{i,j=0}^{N_x}\left[\sum_{k=1}^3(w^n_{k;i,j})^2 \right]\exp(\mu(x_i,y_j)).
	\end{align}
	As initial data we choose $(\tilde{h}, w, v)^T = (1,1,1)^T$ and report on computational results for $L = 3, \Delta x = \Delta y = 1/100$, $L^n, n\geq 0$ for $t_{end} = 3$ and $C_{CFL} = 0.5$.
	\begin{figure}
		\includegraphics[width=0.49\textwidth]{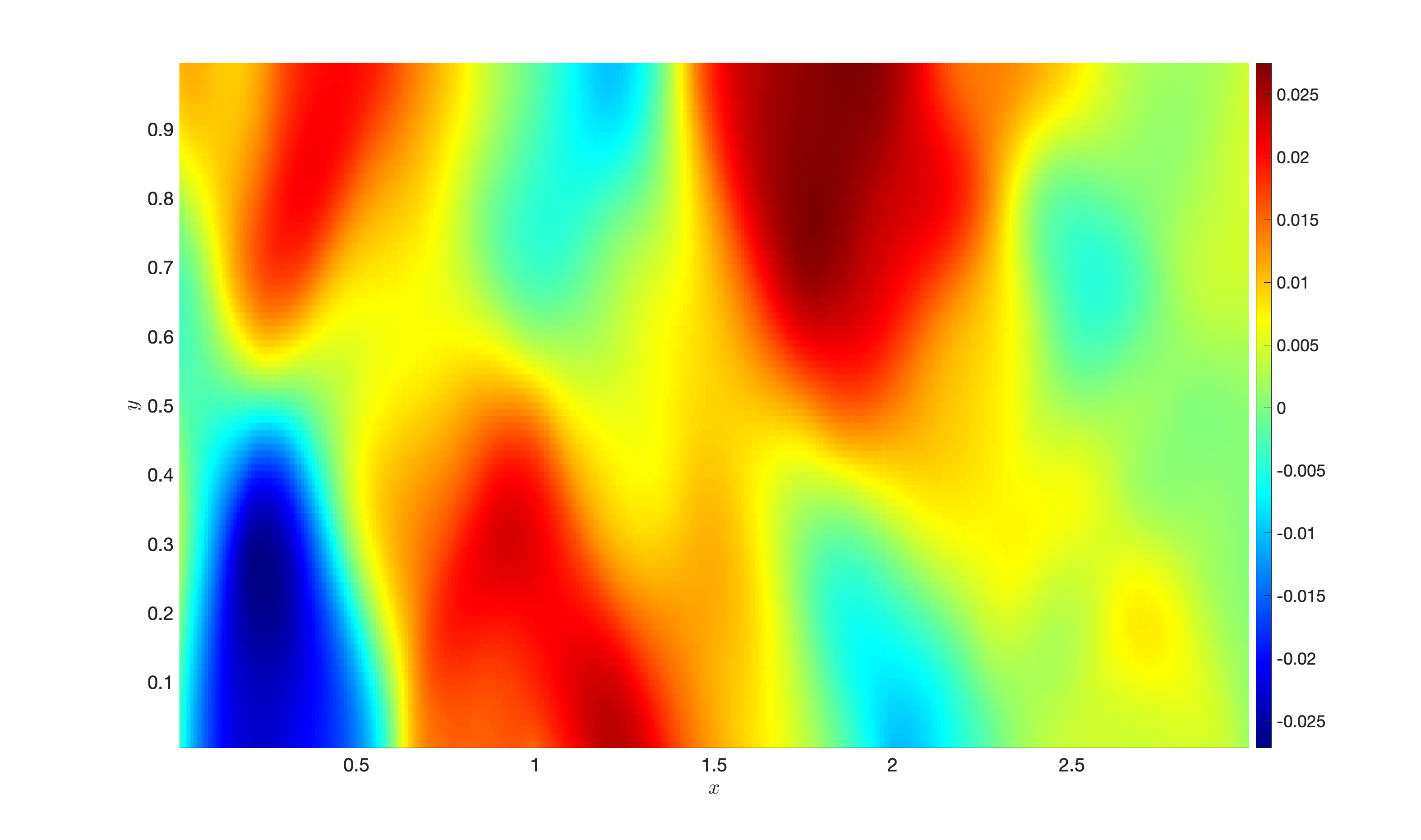}~a)
	    \hfil
	    \includegraphics[width=0.49\textwidth]{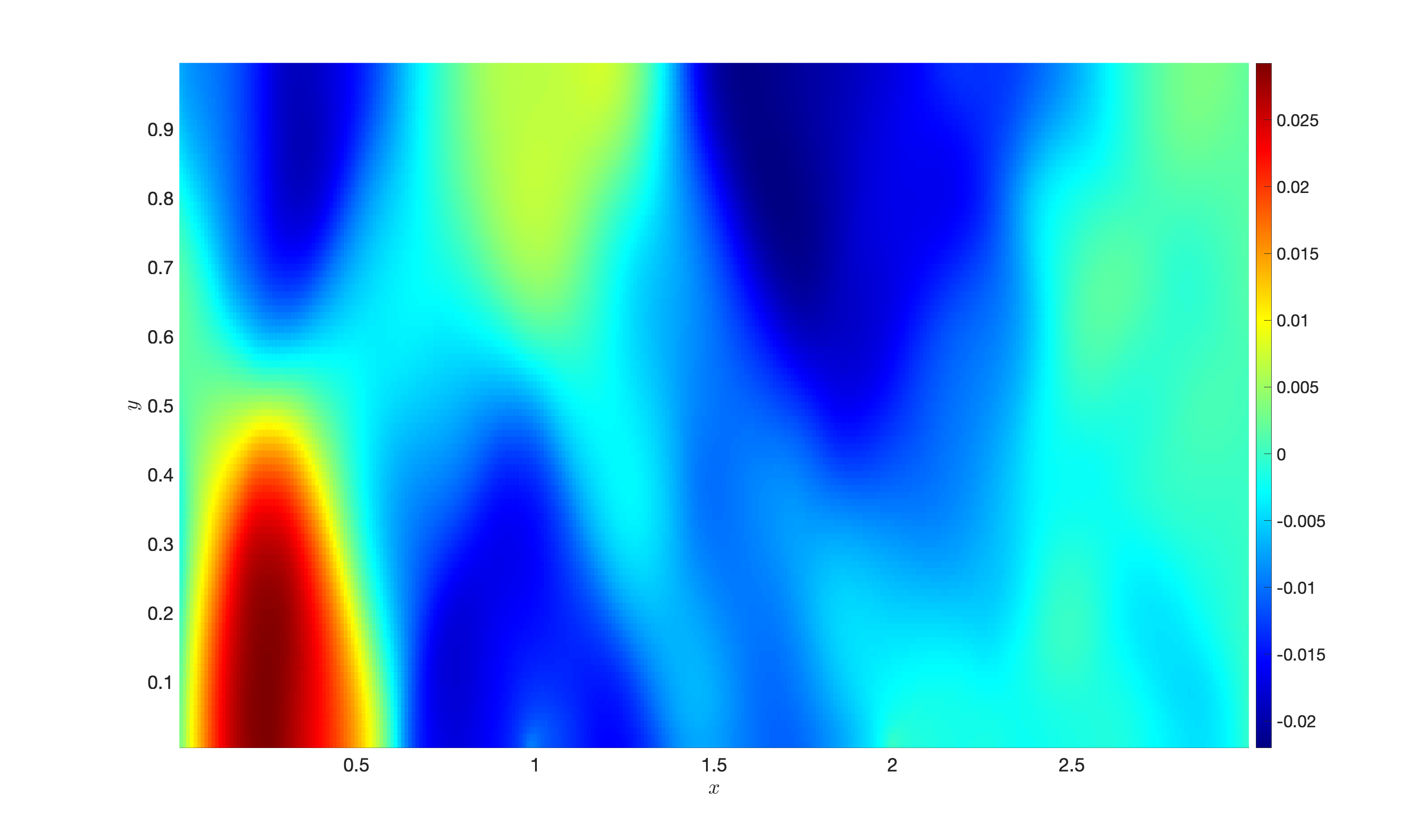}~b)\\
	    \includegraphics[width=0.49\textwidth]{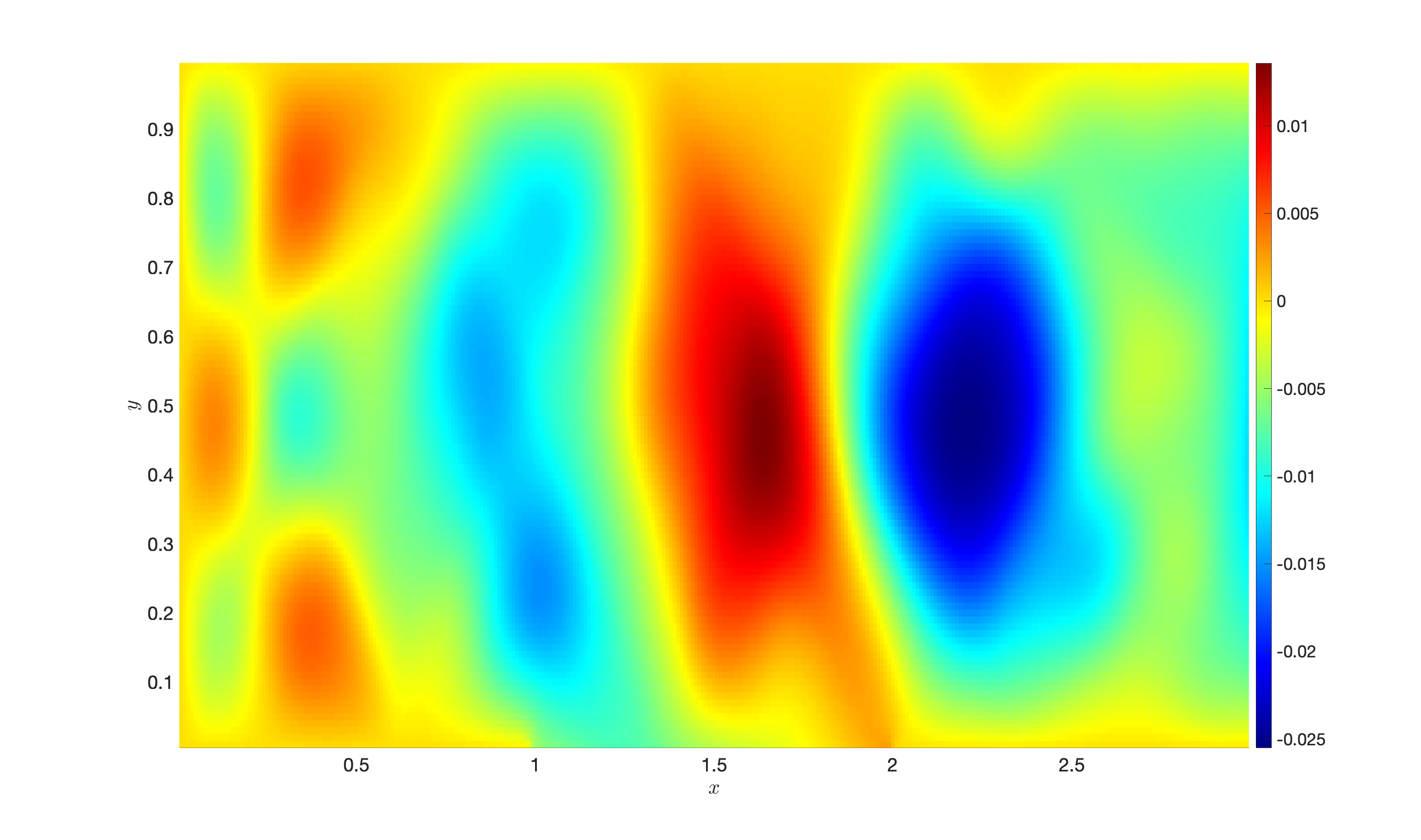}~c)
	    \hfill
	    \includegraphics[width=0.49\textwidth]{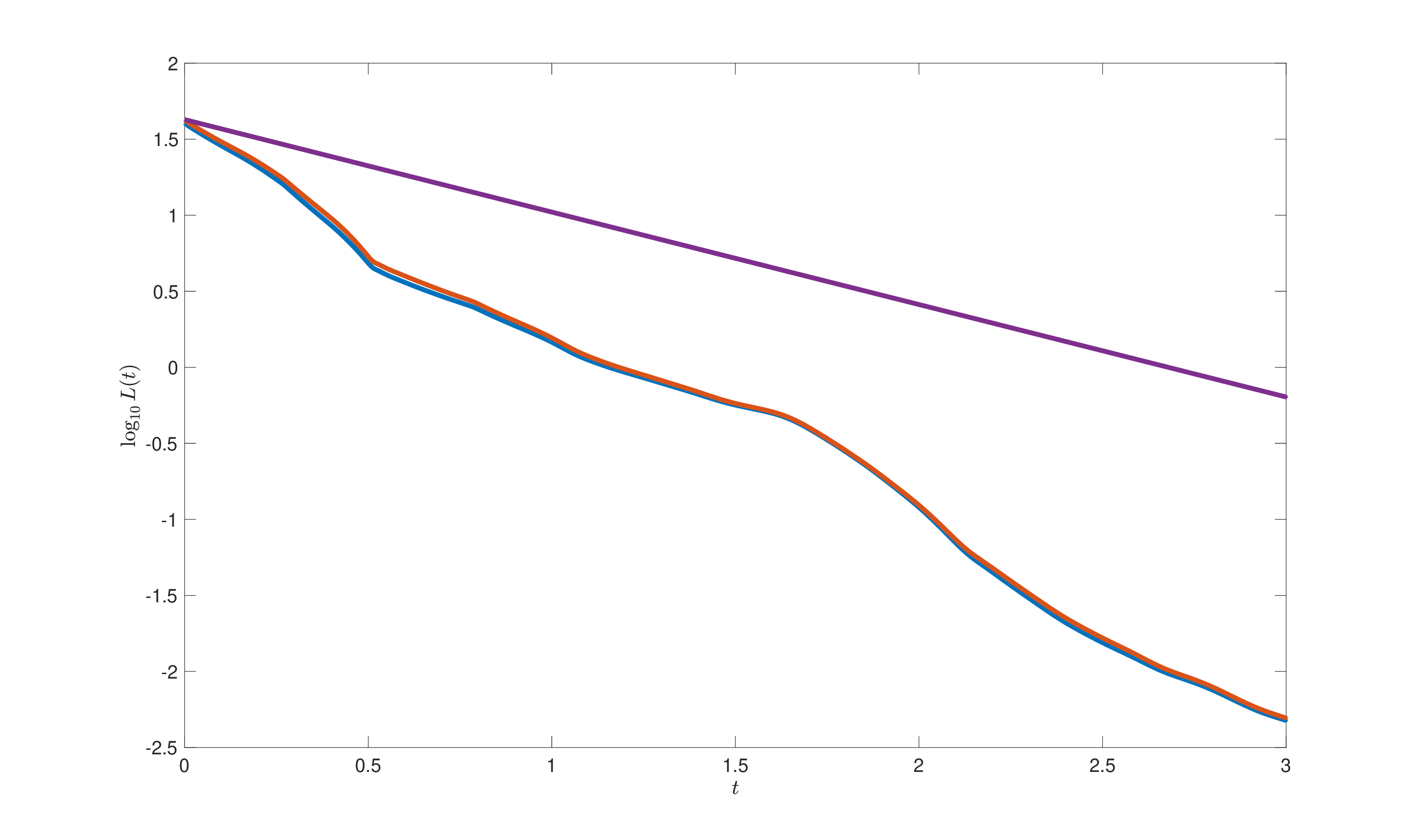}~d)
		\caption{Numerical results at $t_{end} = 3$: The solutions for the components a) $\tilde{h}$, b) $w$ and c) $v$ are given. The computed decay rates for both Lyapunov functions (blue \cite{Yang2024} and red \cite{Herty2023a}) obtained with the MUSCL-FV scheme are given and compared to the decay rate $C = 1.4$ (purple) which is compliant with \eqref{constr_decay_rate}, see d).}
		\label{fig1}
	\end{figure}
	The results are presented in Figure \ref{fig1}.
	It is visible that the observed numerical decay of the Lyapunov function
	 is stronger compared with the theoretical estimate possibly due to additional diffusive terms in the numerical approximation.
	This is confirmed by the observation that coarser grids lead to stronger decay compared with refined meshes, see also \cite{Bambach2022,Thein2023}.
	\subsection{The diagonal system}
	In this section the theoretical estimate on the decay of the Lyapunov function for Example \ref{sec:diag_sys} is confirmed using a numerical discretization of the dynamics.
	The MUSCL second--order finite--volume scheme, see \cite{Toro2009}, is used on a regular mesh for $\Omega$ with grid size $\Delta x = \Delta y$ to solve the discretized equation \eqref{eq:diag_sys_ex}.
	The cell average of $\mathbf{W}$ on $C_{i,j} = [x_i-\frac{\Delta x}2, x_i+\frac{\Delta x}2] \times [y_j-\frac{\Delta x}2, y_j+\frac{\Delta x}2]$, $x_i=i \Delta x, y_j = j\Delta y$ and time $t_n = n\Delta t$ is given by
	\begin{align*}
	    \mathbf{W}_{i,j}^n = \frac{1}{|C_{i,j}|} \int_{C_{i,j}} \mathbf{W}(t_n,x,y)\,\dd x\dd y,
	\end{align*}
	for $i,j=0,\dots,N$ and $n\geq 0$ and where $N \Delta x = 1$. The cell averages of the initial data $\mathbf{W}_0$ are obtained analogously and define $\mathbf{W}_{i,j}^0$.
	As boundary conditions we use transmissive boundary conditions on $\Gamma_1^+, \Gamma^+_2$ and $\Gamma_3^+$, see \cite{Toro2009}.
	Zero boundary conditions are prescribed on $\mathcal{Z}_1, \mathcal{Z}_2$ and $ \mathcal{Z}_3$. The boundary control is imposed on $\mathcal{C}_1, \mathcal{C}_2$ and $\mathcal{C}_3$.
	Here, the control $u^n$ is obtained using a numerical quadrature formula applied to \eqref{eq:u_bc}.
	The Lyapunov function \eqref{eq:lyapunov_gen_diag} is approximated at time $t_n$ by $L^n$ using a numerical quadrature rule on the equi-distant grid 
	\begin{align}
	    L^n:=\Delta x^2 \sum\limits_{i,j=0}^{N_x}\sum_{k=1}^3\left[(w^n_{k;i,j})^2 \exp(\mu_k(x_i,y_j))\right].
	\end{align}
	As initial data a sinusoidal function is chosen, i.e.\ $w_i(0,\mathbf{x}) = \sin(2\pi x)\sin(2\pi y)$.
	We report on computational results for $L^n, n\geq 0$ for the following computational setup $C_L = 4, \Delta x = 10^{-2}, t_{end} = 3$ and $C_{CFL} = 0.5$.
	\begin{figure}
		\includegraphics[width=0.49\textwidth]{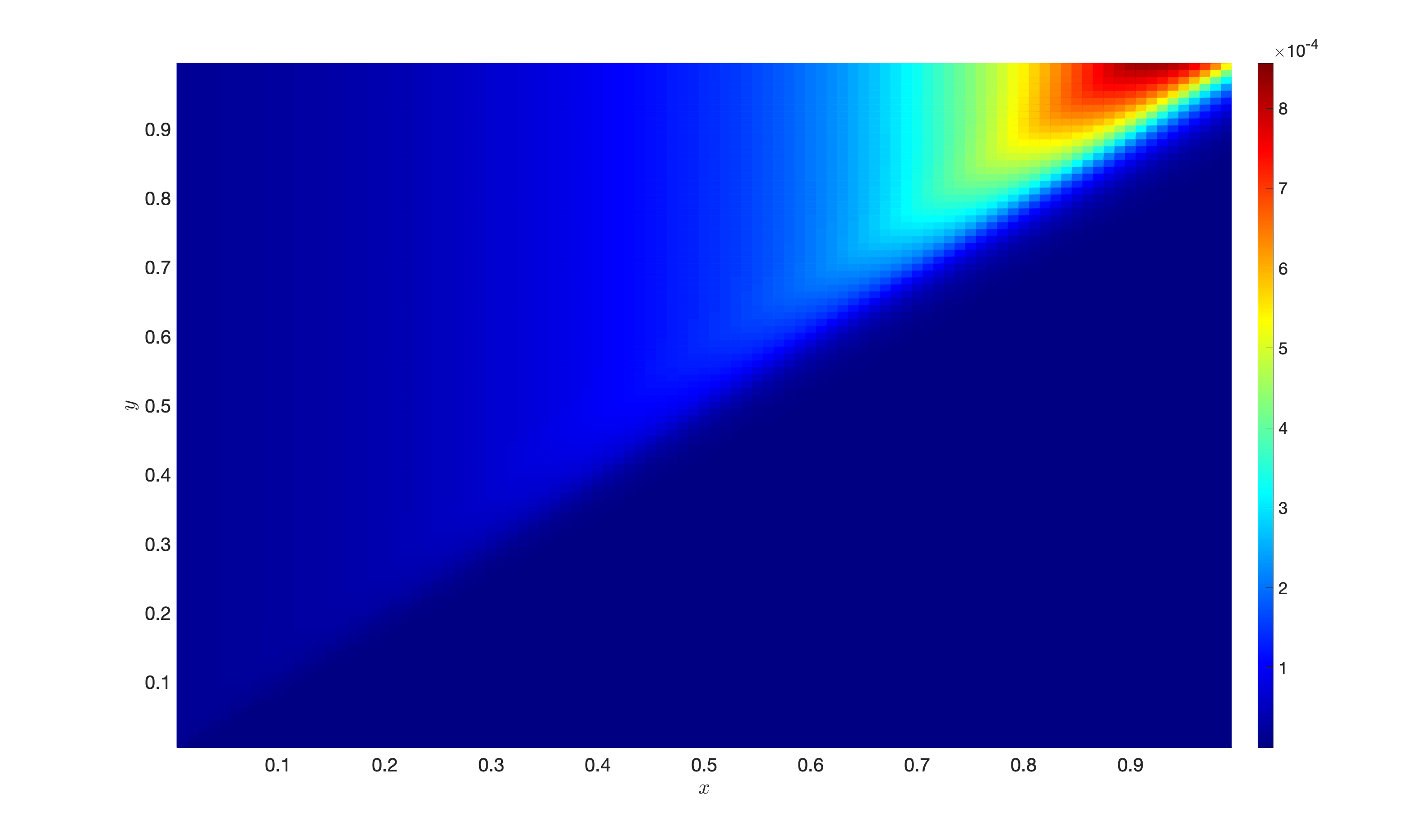}~a)
	    \hfill
	    \includegraphics[width=0.49\textwidth]{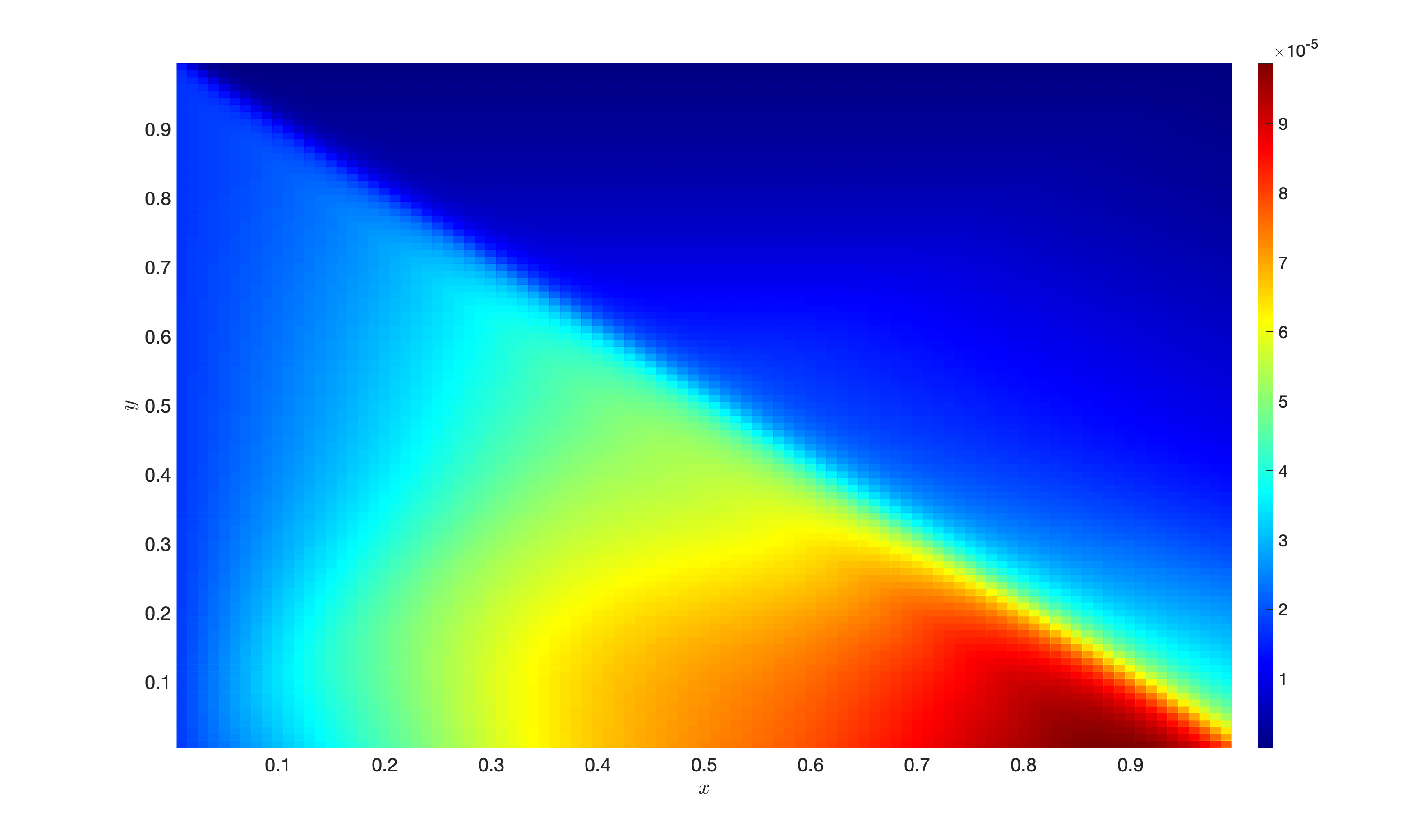}~b)\\
	    \includegraphics[width=0.49\textwidth]{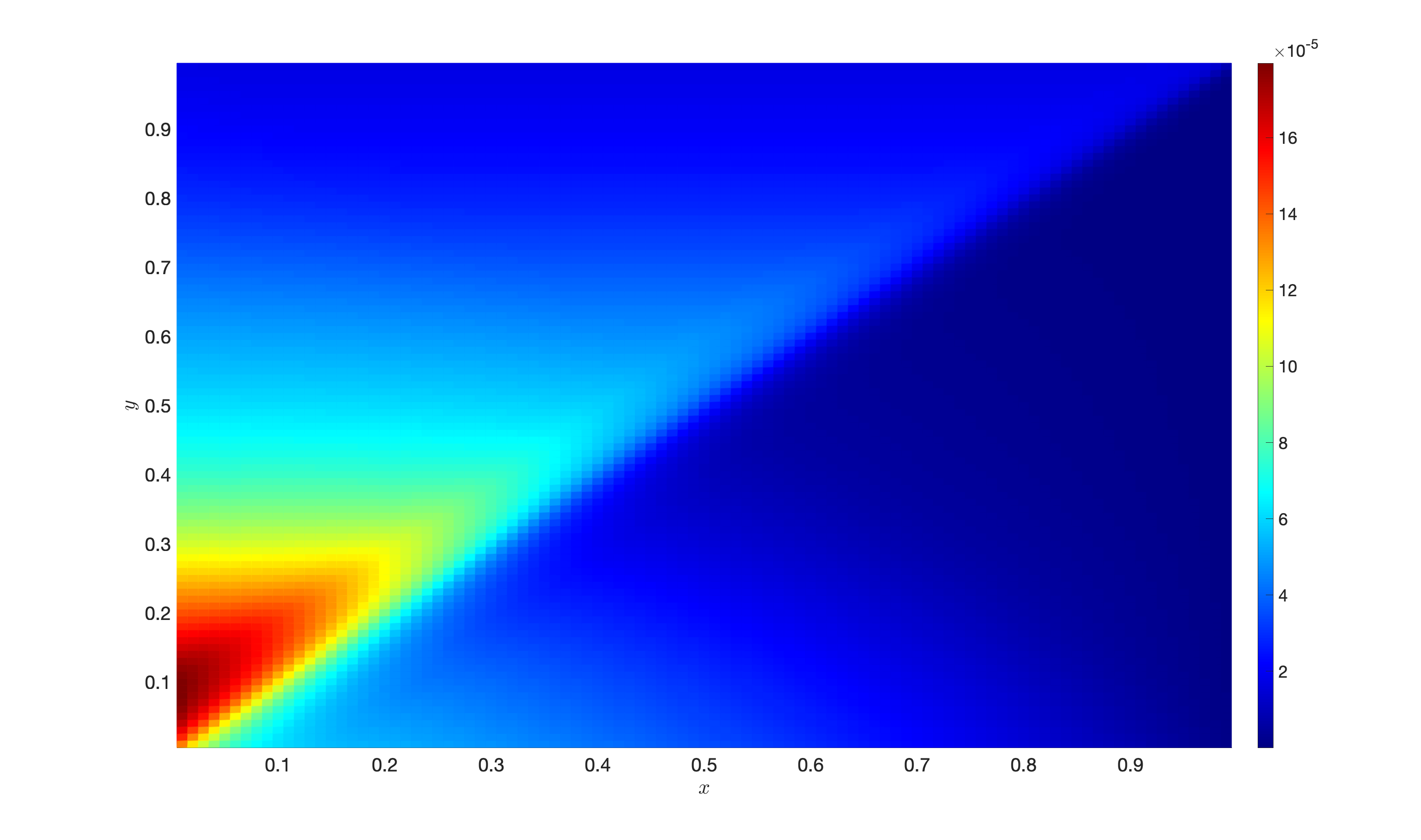}~c)
	    \hfill
	    \includegraphics[width=0.49\textwidth]{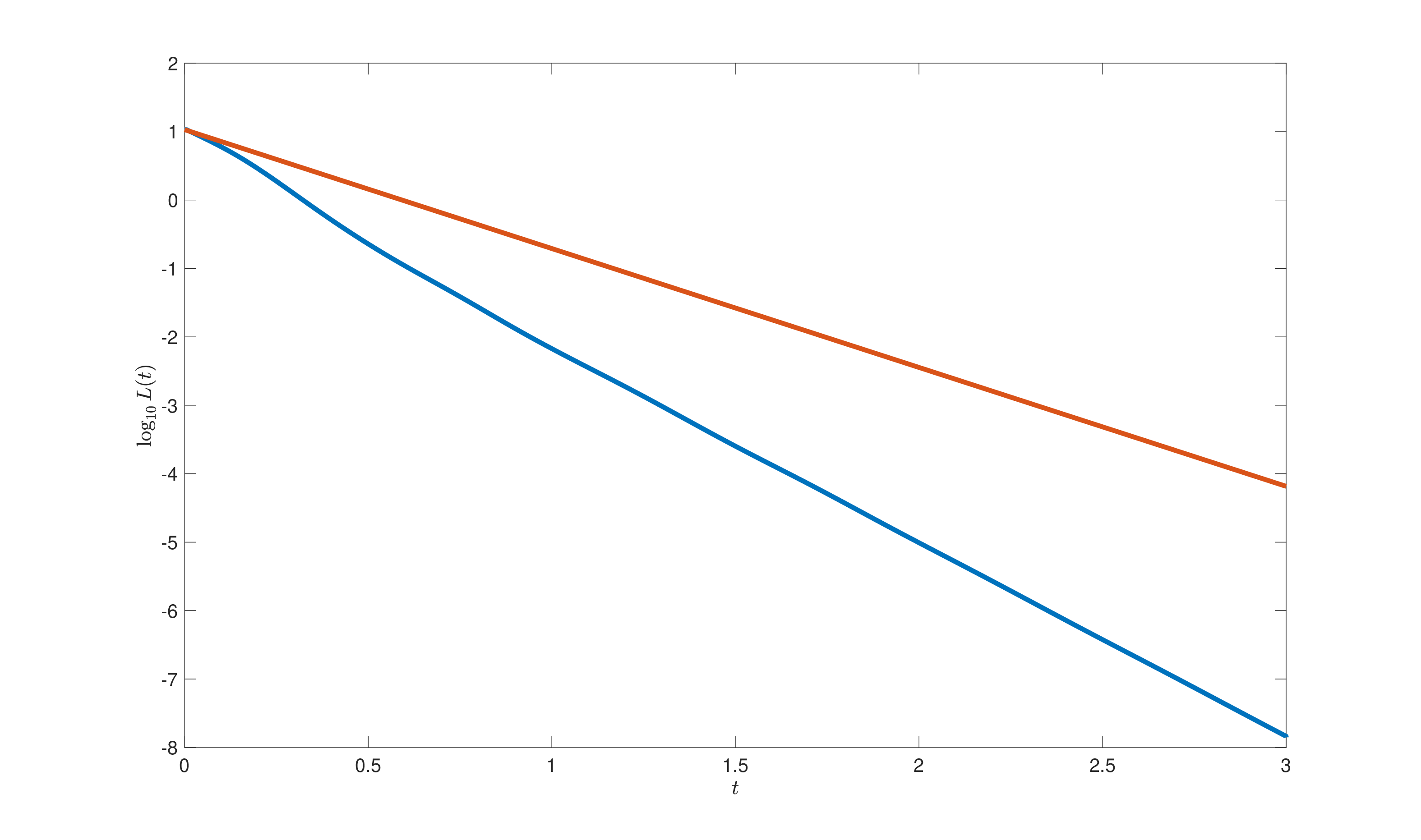}~d)
		\caption{Numerical results at $t_{end} = 3$: The solutions for the components a) $w_1$, b) $w_2$ and c) $w_3$ are given. The computed (blue) decay rate for the Lyapunov function obtained with the MUSCL-FV scheme is given and compared to the theoretical decay rate (red), see d).}
		\label{fig2}
	\end{figure}
	The results are presented in Figure \ref{fig2}.
	It is visible that the observed numerical decay of the Lyapunov function 
	is stronger compared with the theoretical estimate possibly due to additional diffusive terms in the numerical approximation.
	This is confirmed by the observation that coarser grids lead to stronger decay compared with refined meshes, see also \cite{Bambach2022,Thein2023}.
    \section{Summary}
    In the present work we exploit the relation of novel approaches for stabilizing multi--dimensional hyperbolic systems using a boundary feedback control.
    All of which provide exponential decay under suited conditions. The assumptions are briefly reviewed, compared and discussed.
    It is shown that the approach given in \cite{Herty2023a} is always applicable when the prerequisite of \cite{Yang2024} hold.
    Moreover, we highlight the relation between the Lyapunov functions used in these works by showing that a linearization of the weight function used in \cite{Herty2023a} yields the one used in \cite{Yang2024}.
    Finally we discussed the very relevant example of the Saint-Venant equations to highlight our results.
    Moreover, we complete the study by showing that the approach studied in \cite{Herty2023} helps to stabilize systems for which the previous approaches might not be applicable.
    With this work we hope to foster the study of control problems for multi--dimensional hyperbolic systems.
    The present work thus highlights the applicability of the results obtained in \cite{Herty2023a} which is now known to be suitable for the Euler equations, the Saint-Venant equations and systems satisfying the SSC
    condition. These are in particular examples important for different applications.
    It will be interesting to study further complex examples.
    A further point of interest is to determine the Lyapunov potential related to the LMI due to the variety of systems and many possible solutions. It is therefore a subject of further research to identify
    structural properties, besides the SSC given in \cite{Yang2024}, that guarantee the existence of such a function.
    %
    %conflict of interests
    \section*{Conflict of Interest}
    The authors declare that they have no conflict of interest.
    \subsection*{Acknowledgments}
    \small{This research is part of the DFG SPP 2183 \emph{Eigenschaftsgeregelte Umformprozesse}, project 424334423.\\
    M.H. thanks the Deutsche Forschungsgemeinschaft (DFG, German Research Foundation)
    for the financial support through 525842915/SPP2410, 525853336/SPP2410, 320021702/GRK2326, 333849990/IRTG-2379, CRC1481, 423615040/SPP1962, 462234017, 461365406, ERS SFDdM035
    and under Germany’s Excellence Strategy EXC-2023 Internet of Production 390621612 and under the Excellence Strategy of the Federal Government and the L\"ander.\\
    F.T. thanks the Deutsche Forschungsgemeinschaft (DFG, German Research Foundation) for the financial support through 525939417/SPP2410.
    %\\ The authors gratefully acknowledge the valuable comments made by the referee.
    }
    \appendix
    \section{Eigenstructure of the Saint-Venant Equations}\label{app:eigenstruct}
    Subsequently we provide the detailed calculation for the eigenstructure of the considered system \eqref{eq:saintv_sys}.
    The system matrix is given by
    \begin{align}
        \mathbf{A}^\ast(\nu) = \nu_1\mathbf{A}^{(1)} + \nu_2\mathbf{A}^{(2)},\;\nu\in\mathbb{S}^2\label{sv_sys_mat}
    \end{align}
    and thus we have with $v^\ast = 0$
    \[
      \mathbf{A}^\ast(\nu) = \begin{pmatrix} \nu_1w^\ast         & \nu_1\sqrt{gH^\ast} & \nu_2\sqrt{gH^\ast}\\
                                             \nu_1\sqrt{gH^\ast} & \nu_1 w^\ast        & 0\\
                                             \nu_2\sqrt{gH^\ast} & 0                   & \nu_1 w^\ast \end{pmatrix}.
    \]
    The characteristic polynomial is given by
    \begin{align*}
        \chi(\lambda)  &= (\lambda - \nu_1 w^\ast)^3 - (\lambda - \nu_1 w^\ast)\nu_1^2gH^\ast - (\lambda - \nu_1 w^\ast)\nu_2^2gH^\ast\\
        &= (\lambda - \nu_1 w^\ast)\left((\lambda - \nu_1 w^\ast)^2 - \left(\nu_1^2 + \nu_2^2\right)gH^\ast\right) = (\lambda - \nu_1 w^\ast)\left((\lambda - \nu_1 w^\ast)^2 - gH^\ast\right)
    \end{align*}
    and hence we yield the following eigenvalues
    \begin{align}
        \lambda_1(\nu) &= \nu_1 w^\ast - \sqrt{gH^\ast},\;\lambda_2(\nu) = \nu_1 w^\ast,\;\lambda_3(\nu) = \nu_1 w^\ast + \sqrt{gH^\ast}\label{eigvals_sv_sys}\\
        \text{with}\quad\lambda_1(\nu) &< \lambda_2(\nu) < \lambda_3(\nu).\notag
    \end{align}
    The corresponding right eigenvectors are thus obtained to be
    \begin{align}
        \mathbf{R}_1(\nu) = \frac{1}{\sqrt{2}}\begin{pmatrix} 1\\ -\nu_1\\ -\nu_2 \end{pmatrix},\quad
        \mathbf{R}_2(\nu) = \begin{pmatrix} 0\\ -\nu_2\\ \hphantom{-}\nu_1 \end{pmatrix}\quad\text{and}\quad
        \mathbf{R}_3(\nu) = \frac{1}{\sqrt{2}}\begin{pmatrix} 1\\ \nu_1\\ \nu_2 \end{pmatrix}.\label{eigvecs_sv_sys}
    \end{align}
    We hence have the following transformation matrix
    \[
      \mathbf{T}(\nu) = \frac{1}{\sqrt{2}}\begin{pmatrix} 1         & 0                         & 1\\
                                                          -\nu_1    & -\sqrt{2}\nu_2            & \nu_1\\
                                                          -\nu_2    & \hphantom{-}\sqrt{2}\nu_1 & \nu_2 \end{pmatrix}
    \]
    with $\mathbf{\Lambda}^\ast(\nu) = \mathbf{T}^T(\nu)\mathbf{A}^\ast(\nu)\mathbf{T}(\nu)$. The transformation of the state vector $\mathbf{w} = (\tilde{h},w,v)^T$ is given by
    \begin{align}
        \mathbf{v} = \mathbf{T}^T(\nu)\mathbf{w} = \frac{1}{\sqrt{2}}\begin{pmatrix} \tilde{h} - (\nu_1w + \nu_2v)\\ -\sqrt{2}(\nu_2w - \nu_1v)\\ \tilde{h} + \nu_1w + \nu_2v\end{pmatrix}
        \quad\text{and}\quad
        \mathbf{w} = \mathbf{T}(\nu)\mathbf{v} = \frac{1}{\sqrt{2}}\begin{pmatrix} v_1 + v_3\\ -\nu_1(v_1 - v_3) - \sqrt{2}\nu_2 v_2\\ -\nu_2(v_1 - v_3) + \sqrt{2}\nu_1 v_2 \end{pmatrix}
        \label{sv_vec_trafo}
    \end{align}
    \section{Boundary Control}\label{app:bound_contr}
    According to \eqref{bound_split} we can prescribe \textbf{five} controls (counting the boundary parts with negative eigenvalues). In the following we discuss the four boundary parts separately as suggested in \cite{Yang2024}.\\
    \newline
    \textbf{Part I:} We discuss $\{0\}\times[0,1]$ for which we have $\lambda_1 < \lambda_2 < 0 < \lambda_3$ and further
    \[
      \mathbf{v} = \frac{1}{\sqrt{2}}\begin{pmatrix} \tilde{h} + w\\ -\sqrt{2}v\\ \tilde{h} - w\end{pmatrix}.
    \]
    Thus we obtain
    \begin{align}
    	BC^{(I)} = \int_0^1\left[(-w^\ast - \sqrt{gH^\ast})v_1(t,0,y)^2 - w^\ast v_2(t,0,y)^2 + (-w^\ast + \sqrt{gH^\ast})v_3(t,0,y)^2\right]\delta(0,y)\,\dd y \label{boundary_part_1}
    \end{align}
    In view of the present characteristic variables we choose the control
    \begin{align}
    	v_1(t,0,y)^2 = \alpha v_3(t,0,y)^2\quad\text{with}\quad
    	0 \leq \alpha \leq \frac{\sqrt{gH^\ast} - w^\ast}{\sqrt{gH^\ast} + w^\ast} < 1.\label{contr:boundary_part_1}
    \end{align}
    With this choice we can conclude
    \begin{align}
    	BC^{(I)} \geq -w^\ast\int_0^1 v_2(t,0,y)^2 \delta(0,y)\,\dd y \label{reminder:boundary_part_1}
    \end{align}
    and we will treat this term later on.\\
    \newline
    \textbf{Part II:} We discuss $[0,L]\times\{0\}$ for which we have $\lambda_1 < \lambda_2 = 0 < \lambda_3$ and further
    \[
      \mathbf{v} = \frac{1}{\sqrt{2}}\begin{pmatrix} \tilde{h} + v\\ \sqrt{2}w\\ \tilde{h} - v\end{pmatrix}.
    \]
    Thus we obtain
    \begin{align}
    	BC^{(II)} = \int_0^L\left[(-\sqrt{gH^\ast})v_1(t,x,0)^2 + \sqrt{gH^\ast}v_3(t,x,0)^2\right]\delta(x,0)\,\dd x \label{boundary_part_2}
    \end{align}
    The spillway is located at $x \in [L/3,2L/3]$ and the boundary part $[0,L/3)\cup (2L/3,L]$ is assumed to be a solid wall. For the solid wall we thus assume the normal velocity to be zero, i.e. for all $x \in [0,L/3)\cup (2L/3,L]$ we have the following control
    \begin{align}
    	v(t,x,0) = 0\quad\Leftrightarrow\quad (\tilde{h} + v)(t,x,0) = (\tilde{h} - v)(t,x,0)\quad\Rightarrow\quad v_1(t,x,0)^2 = v_3(t,x,0)^2.\label{contr:boundary_part_2}
    \end{align}
    With this choice we can conclude
    \begin{align}
    	BC^{(II)} = \sqrt{gH^\ast}\int_{L/3}^{2L/3} (v_3(t,x,0)^2 - v_1(t,x,0)^2) \delta(x,0)\,\dd x \label{reminder:boundary_part_2}
    \end{align}
    and we will treat this term later on.\\
    \newline
    \textbf{Part III:} We discuss $\{L\}\times[0,1]$ for which we have $\lambda_1 < 0 < \lambda_2 < \lambda_3$ and further
    \[
      \mathbf{v} = \frac{1}{\sqrt{2}}\begin{pmatrix} \tilde{h} - w\\ \sqrt{2}v\\ \tilde{h} + w\end{pmatrix}.
    \]
    Thus we obtain
    \begin{align}
    	BC^{(III)} = \int_0^1\left[(w^\ast - \sqrt{gH^\ast})v_1(t,L,y)^2 + w^\ast v_2(t,L,y)^2 + (w^\ast + \sqrt{gH^\ast})v_3(t,L,y)^2\right]\delta(L,y)\,\dd y \label{boundary_part_3}
    \end{align}
    In view of the present characteristic variables we choose the control
    \begin{align}
    	v_1(t,L,y)^2 = \beta v_3(t,L,y)^2\quad\text{with}\quad
    	0 \leq \beta \leq \frac{\sqrt{gH^\ast} + w^\ast}{\sqrt{gH^\ast} - w^\ast}. \label{contr:boundary_part_3}
    \end{align}
    With this choice we can conclude
    \begin{align*}
    	BC^{(III)} \geq w^\ast\int_0^1 v_2(t,L,y)^2 \delta(L,y)\,\dd y \geq 0.
    \end{align*}
    \textbf{Part IV:} We discuss $[0,L]\times\{1\}$ for which we have $\lambda_1 < \lambda_2 = 0 < \lambda_3$ and further
    \[
      \mathbf{v} = \frac{1}{\sqrt{2}}\begin{pmatrix} \tilde{h} - v\\ -\sqrt{2}w\\ \tilde{h} + v\end{pmatrix}.
    \]
    Thus we obtain
    \begin{align}
    	BC^{(IV)} = \int_0^L\left[(-\sqrt{gH^\ast})v_1(t,x,1)^2 + \sqrt{gH^\ast}v_3(t,x,1)^2\right]\delta(x,1)\,\dd x \label{boundary_part_4}
    \end{align}
    The boundary part (IV) is assumed to be a solid wall and we thus assume the normal velocity to be zero, i.e. for all $x \in [0,L]\times\{1\}$ we have the following control
    \begin{align}
    	v(t,x,1) = 0\quad\Leftrightarrow\quad (\tilde{h} - v)(t,x,1) = (\tilde{h} + v)(t,x,1)\quad\Rightarrow\quad v_1(t,x,1)^2 = v_3(t,x,1)^2. \label{contr:boundary_part_4}
    \end{align}
    With this choice we can conclude
    \begin{align*}
    	BC^{(IV)} = 0.
    \end{align*}
    Summing up we yield with \eqref{reminder:boundary_part_1} and \eqref{reminder:boundary_part_2}
    \begin{align*}
    	\mathcal{BC} &= BC^{(I)} + BC^{(II)} + \underbrace{BC^{(III)}}_{\geq 0} + \underbrace{BC^{(IV)}}_{=0}\\
    	&\geq -w^\ast\int_0^1 v_2(t,0,y)^2 \delta(0,y)\,\dd y + \sqrt{gH^\ast}\int_{L/3}^{2L/3} (v_3(t,x,0)^2 - v_1(t,x,0)^2) \delta(x,0)\,\dd x\\
    	&= -w^\ast\int_0^1 v(t,0,y)^2 \delta(0,y)\,\dd y + \sqrt{gH^\ast}\int_{L/3}^{2L/3} \left(\frac{1}{2}[\tilde{h} - v]^2(t,x,0) - \frac{1}{2}[\tilde{h} + v]^2(t,x,0)\right) \delta(x,0)\,\dd x\\
    	&=: \mathcal{R}.
    \end{align*}
    We need to choose the controls for $v(t,0,y)$ and $[\tilde{h} + v](t,x,0)$ such that $\mathcal{R}$ becomes non-negative. Several approaches are now possible. We follow \cite{Yang2024} and the literature cited therein and want to control the normal velocity depending the measured water height and thus use $\tilde{h}$. Let us rewrite $\mathcal{R}$ as follows using $x = L(y + 1)/3$
    \begin{align*}
    	\mathcal{R} &= -w^\ast\int_0^1 v^2(t,0,y) \delta(0,y)\,\dd y\\
    	&+ \frac{L}{3}\sqrt{gH^\ast}\int_0^1 \left(\frac{1}{2}[\tilde{h} - v]^2\left(t,\frac{L}{3}(y + 1),0\right) - \frac{1}{2}[\tilde{h} + v]^2\left(t,\frac{L}{3}(y + 1),0\right)\right) \delta\left(\frac{L}{3}(y + 1),0\right)\,\dd y
    \end{align*}
    We now prescribe the following controls
    \begin{align}
    	v_2^2(t,0,y) &= \gamma v_3^2\left(t,\frac{L}{3}(y+1),0\right)\quad\Leftrightarrow\quad v^2(t,0,y) = \gamma \frac{1}{2}[\tilde{h} - v]^2\left(t,\frac{L}{3}(y+1),0\right) \label{contr_b:boundary_part_1}\\
    	v_1^2(t,x,0) &= \varepsilon v_3^2(t,x,0)\quad\Leftrightarrow\quad [\tilde{h} + v]^2(t,x,0) = \varepsilon [\tilde{h} - v]^2(t,x,0)\label{contr_b:boundary_part_2}\\
    	\text{with}\quad \varepsilon &\in [0,1]\quad\text{and}\quad 0 \leq \gamma \leq (1 - \varepsilon)\frac{2L}{9}\frac{\sqrt{gH^\ast}}{w^\ast} \label{control_conditions}    	
    \end{align}
    The condition on $\gamma$ and $\varepsilon$ ensures the following in the case of \cite{Herty2023a} and \cite{Yang2024}
    \begin{align*}
    	&\phantom{=}-\gamma w^\ast\delta(0,y) + \frac{L}{3}(1 - \varepsilon)\sqrt{gH^\ast}\delta\left(\frac{L}{3}(y + 1),0\right)\\
    	&\geq -2L\gamma w^\ast + \frac{L}{3}(1 - \varepsilon)\sqrt{gH^\ast}\left(2L - \frac{L}{3}(y + 1)\right)\\
    	&\geq -2L\gamma w^\ast + \frac{L}{3}(1 - \varepsilon)\sqrt{gH^\ast}\left(2L - \frac{2L}{3}\right)\\
    	&= -2L\gamma w^\ast + \frac{4L^2}{9}(1 - \varepsilon)\sqrt{gH^\ast}\\
    	&\geq 0.
    \end{align*}
    Thus we have the positivity of the boundary integral. Note that for the positivity only the relation between the squares of the variables is needed and there is some freedom when prescribing the controls for the variables themselves. We therefore follow \cite{Yang2024}. The basic idea is to adjust the water velocity according to the height such that the height is increased via water inflow and decreased via water outflow, respectively.
    We want to summarize the five controls in terms of the characteristic variables.
    \begin{enumerate}[(I)]
    	\item $\{0\}\times[0,1]$: The controls are determined according to \eqref{contr:boundary_part_1}, \eqref{contr_b:boundary_part_1} and \eqref{control_conditions}
    	\begin{align}
    		v_1(t,0,y) &= \sqrt{\alpha}v_3(t,0,y),\label{contr_a_charvar_bp1}\\
    		v_2(t,0,y) &= -\sqrt{\gamma}v_3\left(t,\frac{L}{3}(y + 1),0\right)\label{contr_b_charvar_bp1}.
    	\end{align}
    	\item $[0,L]\times\{0\}$: The controls are determined according to \eqref{contr:boundary_part_2}, \eqref{contr_b:boundary_part_2} and \eqref{control_conditions}
    	\begin{align}
    		v_1(t,x,0) &= v_3(t,x,0), x \in [0,L/3)\cup(2L/3,L] \label{contr_a_charvar_bp2}\\
    		v_1(t,x,0) &= \sqrt{\varepsilon}v_3(t,x,0), x \in [L/3,2L/3].\label{contr_b_charvar_bp2}
    	\end{align}
    	\item $\{1\}\times[0,1]$: The controls are determined according to \eqref{contr:boundary_part_3}
    	\begin{align}
    		v_1(t,L,y) &= \sqrt{\beta}v_3(t,L,y).\label{contr_charvar_bp3}
    	\end{align}
		\item $[0,L]\times\{1\}$: The controls are determined according to \eqref{contr:boundary_part_4}
    	\begin{align}
    		v_1(t,x,1) &= v_3(t,x,1).\label{contr_charvar_bp4}
    	\end{align}
    \end{enumerate}
    To obtain the controls in the balance (primitive) variables we first use \eqref{contr_b_charvar_bp2} to obtain
    \begin{align*}
    	v_1(t,x,0) = \sqrt{\varepsilon}v_3(t,x,0)\quad\Leftrightarrow\quad v(t,x,0) = -\frac{1 - \sqrt{\varepsilon}}{1 + \sqrt{\varepsilon}}\tilde{h}(t,x,0).
    \end{align*}
    This can now be inserted in \eqref{contr_b_charvar_bp1} to obtain
    \begin{align*}
    	v_2(t,0,y) &= -\sqrt{\gamma}v_3\left(t,\frac{L}{3}(y + 1),0\right)\\
    	\Leftrightarrow\quad v(t,0,y) &= \sqrt{\frac{\gamma}{2}}[\tilde{h} - v]\left(t,\frac{L}{3}(y + 1),0\right)\\
    	\Leftrightarrow\quad v(t,0,y) &= \sqrt{\frac{\gamma}{2}}\left[1 + \frac{1 - \sqrt{\varepsilon}}{1 + \sqrt{\varepsilon}}\right]\tilde{h}\left(t,\frac{L}{3}(y + 1),0\right)\\
    	\Leftrightarrow\quad v(t,0,y) &= \frac{\sqrt{2\gamma}}{1 + \sqrt{\varepsilon}}\tilde{h}\left(t,\frac{L}{3}(y + 1),0\right).
    \end{align*}
    We hence obtain the five controls in terms of the primitive variables
    \begin{enumerate}[(I)]
    	\item $\{0\}\times[0,1]$: The controls are determined according to \eqref{contr:boundary_part_1}, \eqref{contr_b:boundary_part_1} and \eqref{control_conditions}
    	\begin{align}
    		w(t,0,y) &= -\frac{1 - \sqrt{\alpha}}{1 + \sqrt{\alpha}}\tilde{h}(t,0,y),\label{contr_a_primvar_bp1}\\
    		v(t,0,y) &= \frac{\sqrt{2\gamma}}{1 + \sqrt{\varepsilon}}\tilde{h}\left(t,\frac{L}{3}(y + 1),0\right)\label{contr_b_primvar_bp1}.
    	\end{align}
    	\item $[0,L]\times\{0\}$: The controls are determined according to \eqref{contr:boundary_part_2}, \eqref{contr_b:boundary_part_2} and \eqref{control_conditions}
    	\begin{align}
    		v(t,x,0) &= 0, x \in [0,L/3)\cup(2L/3,L] \label{contr_a_charvar_bp2}\\
    		v(t,x,0) &= -\frac{1 - \sqrt{\varepsilon}}{1 + \sqrt{\varepsilon}}\tilde{h}(t,x,0), x \in [L/3,2L/3].\label{contr_b_primvar_bp2}
    	\end{align}
    	\item $\{1\}\times[0,1]$: The controls are determined according to \eqref{contr:boundary_part_3}
    	\begin{align}
    		w(t,L,y) &= \frac{1 - \sqrt{\beta}}{1 + \sqrt{\beta}}\tilde{h}(t,L,y).\label{contr_primvar_bp3}
    	\end{align}
		\item $[0,L]\times\{1\}$: The controls are determined according to \eqref{contr:boundary_part_4}
    	\begin{align}
    		v(t,L,y) &= 0.\label{contr_primvar_bp4}
    	\end{align}
    \end{enumerate}
    Concerning \cite{Yang2024} we thus have the same controls setting
    \[
      k_1 = \sqrt{\beta},\;k_2 = \sqrt{\alpha},\;k_3 = \sqrt{\varepsilon}\;\text{and}\;k_4 = \sqrt{\frac{\gamma}{2}}
    \]
    % Literature
    \phantomsection
    \bibliographystyle{abbrvurl}
    \bibliography{comp_bfc_lit.bib}
\end{document}